\documentclass[fleqn]{article}
 \usepackage{amscd}
                     \usepackage{amssymb}
 \newcommand{\be}{\begin{equation}}
       \newcommand{\ee}{\end{equation}}
       \newcommand{\ba}{\begin{eqnarray}}
        \newcommand{\ea}{\end{eqnarray}}
 \newcommand{\ban}{\begin{eqnarray*}}
 \newcommand{\ean}{\end{eqnarray*}}

 \newcommand{\lp}{\langle}
 \newcommand{\rp}{\rangle}
 \newcommand{\ra}{\rightarrow}

  \newcommand{\qed}{\hspace*{\fill}\rule{3mm}{3mm}\quad \vspace{.2cm}}
  \newcommand{\Pf}{\noindent {\bf Proof:} }
  \newcommand{\Rk}{\noindent {\bf Remark} }

 \newcommand{\sect}[1]{\section{#1} \setcounter{equation}{0}}
 
 \newcommand{\gap}{\mbox{gap}}

 \newtheorem{theo}{Theorem}[section]
 
\begin{document}
 \newtheorem{defn}[theo]{Definition}
 \newtheorem{ques}[theo]{Question}
 \newtheorem{lem}[theo]{Lemma}
 \newtheorem{lemma}[theo]{Lemma}
 \newtheorem{prop}[theo]{Proposition}
 \newtheorem{coro}[theo]{Corollary}
 \newtheorem{ex}[theo]{Example}
 \newtheorem{note}[theo]{Note}
 \newtheorem{conj}[theo]{Conjecture}
 
 \title{The covering spectrum of a compact length space}
 \author{Christina Sormani\thanks {Partially supported by NSF
 Grant \# DMS-0102279 and a grant from the City University
 of New York PSC-CUNY Research Award Program}   \and Guofang Wei\thanks
 {Partially
 supported by NSF Grant \# DMS-0204187.}} \date{}
 \maketitle

 \begin{abstract}
We define a new spectrum for compact length spaces and
Riemannian manifolds called the ``covering spectrum"
which roughly measures the size of the one dimensional
holes in the space.   More specifically, the covering
spectrum is a set of real numbers $\delta>0$ 
which identify the distinct $\delta$ covers of the
space.  We investigate the relationship between this
covering spectrum, the length spectrum, the marked length
spectrum and the Laplace spectrum.  We analyze the
behavior of the covering spectrum under Gromov-Hausdorff 
convergence and study its gap phenomenon.
  \end{abstract}
 
 \newcommand{\inj}{\mbox{inj}}
 \newcommand{\vol}{\mbox{vol}}
 \newcommand{\diam}{\mbox{diam}}
 \newcommand{\Ric}{\mbox{Ric}}
 \newcommand{\Iso}{\mbox{Iso}}
 \newcommand{\Hess}{\mbox{Hess}} 
 \newcommand{\divg}{\mbox{div}}
 
 \sect{Introduction}

One of the most important subfields of Riemannian Geometry is the study of 
the Laplace spectrum of a compact Riemannian manifold.  Recall that the 
Laplace spectrum is defined as the set of eigenvalues of the Laplace 
operator.  The elements of the Laplace spectrum are assigned a multiplicity
equal to the dimension of the corresponding eigenspace.

Another spectrum defined in an entirely 
different manner is the length spectrum of a manifold: the set of 
lengths of smoothly closed geodesics.  There are various methods used to
assign a multiplicity to each element of the length 
spectrum.  The simplest notion
is to count all geodesics of a given length. This becomes uninteresting
when one has continua of geodesics of the same length as in a torus, so that 
all or some
multiplicities become infinite. A common alternative definition of the
multiplicity of a given length
is the number 
of free homotopy classes of geodesics that contain a smoothly closed 
geodesic sharing that length (c.f. \cite{Go1}).  We will use the latter
definition.

It was proven by Colin de Verdiere \cite{CdV} that the Laplace spectrum
determines the length spectrum of a generic manifold.  
(See also Duistermaat-Guillemin's paper \cite{DuGu}).  
In particular, the
Laplace spectrum determines the length spectrum on negatively curved manifolds 
of arbitrary dimension.  
However, there are pairs of isospectral manifolds first
constructed by Carolyn Gordon \cite{Go1} that have different length spectra 
when one takes multiplicity into account.  These pairs are Heisenberg manifolds
and
Pesce has since shown that the length spectrum, not counting multiplicity,
is determined by the Laplace spectrum on Heisenberg manifolds \cite{Ps}. 

There is also a concept called the
marked length spectrum which gives the 
lengths of smoothly closed geodesics freely homotopic to a representative
of each element in the fundamental group. 
One has the remarkable result that compact surfaces of negative curvature
with same marked length spectrum are isometric \cite{Ot, Cr, Fa}.  This is not
true in general, as the sphere and the Zoll
sphere have same marked length spectrum but are not isometric \cite{Bes}.
Gornet has shown that 
Laplace isospectral nilmanifolds with the same marked length spectrum need
not be isometric or have the same spectrum on one forms \cite{Grnt}. 

In this paper, we have defined a new spectrum for compact Riemannian 
manifolds which we call the {\em covering spectrum}  [See
Defn~\ref{coveringspectrum}].  
In fact, this spectrum can be well defined on compact length spaces 
[Defn~\ref{compactlength}].  Note that it isn't too difficult to extend the
concept of a length spectrum to such spaces, but there is no natural Laplace 
spectrum unless one adds an appropriately defined measure on the metric space,
(see e.g. \cite[Section 6]{ChCo3}). 

The authors first defined a special sequence of covering spaces 
for a given complete length space, $X$, called the delta covers of $X$
in \cite{SoWei} (see Defn~\ref{defdel}).  We used 
these delta covers to study the fundamental groups of these spaces and 
their universal covers.  In particular, we proved that the universal cover 
of a compact length space $X$ is a $\delta$ cover for a sufficiently small 
real number $\delta$ \cite[Proposition 3.2]{SoWei}.  We can now show that a
compact length $X$ has a universal cover iff there is only a finite set of
distinct delta covers [Thm~\ref{univ}].  We have named the 
corresponding finite list of distinct real numbers the ``covering spectrum" 
of $X$ [Defn~\ref{coveringspectrum}].  Roughly, this covering spectrum lists 
the sizes of one 
dimensional holes in $X$.  For example, the covering spectrum of a $1\times3$
flat torus is $\{1/2, 3/2\}$ and the covering spectrum of the standard 
$\mathbb RP^2$ is $\{\pi/2\}$.  Recall 
that if $X$ is a compact Riemannian manifold, then it is a length space and 
has a universal cover, so its covering spectrum is well defined and finite.  

Compact length spaces have grown in interest among Riemannian Geometers 
in recent years because they are the natural limits of Riemannian 
manifolds using Gromov's compactness theorem \cite{Gr}.  Gromov-Hausdorff
limits of Cauchy sequences of Riemannian manifolds with a uniform upper bound on
diameter are compact length 
spaces and, with appropriate curvature bounds on the manifolds, they are 
metric measure spaces \cite{Gr} \cite{ChCo1}.  Cheeger and Colding have 
proven Fukaya's conjecture
that the Laplace spectra of a sequence of manifolds with a uniform lower Ricci
curvature 
bound converge to the Laplace spectrum of the metric measure limit space
\cite{ChCo3}.  It is important to note that one needs metric measure
convergence of the manifolds not just Gromov-Hausdorff convergence
to control the Laplace spectrum in this way \cite{Fu}.

On the other hand Gromov-Hausdorff convergence does not interact well with 
length spectra in general.  This is because closed geodesics can disappear and
appear in the limit and the length spectrum of the sequence doesn't converge 
to the length spectrum of the limit 
(c.f. Examples~\ref{closed-geodesics1}-~\ref{closed-geodesics3}).  

Here we have shown that the  covering spectrum interacts very nicely 
with Gromov-Hausdorff convergence [Thm~\ref{covspecconv}] and is 
fairly easy to define both on 
manifolds and limit spaces.  This follows from the fact that the delta
covering spaces are well controlled when the base spaces converge in the
Gromov-Hausdorff sense
[Theorem~\ref{delta-precom}] and \cite[Theorem 3.6]{SoWei}. 
Another interesting property is that the 
covering spectrum, when assigned an appropriate multiplicity, may be 
used to study fundamental groups [Defn~\ref{multi}, Prop~\ref{mcpg}].

We prove that every element in the covering spectrum is (1/2) of an element in
the length spectrum [Thm~\ref{cov-length-spectrum}].  We also prove that the
marked length
spectrum determines the covering spectrum  
on any compact length space with a universal cover [Thm~\ref{mlscs}] .
We also discuss the relationship between  the covering and the Laplace spectra 
on compact Riemannian manifolds and give a number of examples as described
below.
In particular we construct Laplace isospectral Heisenberg manifolds with
different
covering spectra [Ex~\ref{isolapnotcov}].

The paper is organized as follows:

In Section~\ref{background} we provide all the necessary background
including the definition of delta covers and some key examples.
In particular we recall that a universal cover is a cover of all covers
and that the Hawaii ring is a compact length space with no universal cover.

In Section~\ref{delcovsect} we define the covering spectrum,
$CovSpec(X)$, for an arbitrary compact length space, $X$, and prove that
$CovSpec(X)$ is discrete and $Cl(CovSpec(X) \subset \mathbb R)\subset
CovSpec(X)\cup\{0\}$ [Prop~\ref{discrete}].
We then prove that $CovSpec(X)$ is finite iff $X$ has a universal cover
[Thm~\ref{univ}].

In Section~\ref{covandlength} we restrict our attention to compact length 
spaces
that have a universal cover.  We extend the definition of length spectrum,
[Defn~\ref{lsdef}to these spaces and 
prove Theorem~\ref{cov-length-spectrum} 
that $CovSpec(X) \subset (1/2) LengthSpec(X)$.
We then further restrict ourselves to compact length spaces with
simply connected universal covers and extend the definition
of the minimal length spectrum [Defn~\ref{mlsdef}].
We prove this spectrum is closed and discrete
and that $CovSpec(X) \subset (1/2) MinLengthSpec(X)$ 
[Theorem~\ref{2cov-length-spectrum}] .

In Section~\ref{markcovandlen} we extend the definition
of the marked length spectrum to these compact length
spaces with universal covers [Defn~\ref{markdef}].   
Note that to extend the definition of the
marked length spectrum which ordinarily depends on the fundamental group of
the manifold, we use the ``revised fundamental group'' instead.  
This is the group
of deck transforms
of the universal cover [Defn~\ref{revfg}].  

We then prove Theorem~\ref{mlscs} that the marked length 
spectrum determines the covering
spectrum.  In fact, Thm~\ref{mlscs} also relates the 
covering spectrum to a special
sequence of subgroups of the revised fundamental group, 
which is then used
to define multiplicity for the covering spectrum in Section~\ref{sectmult}.  

As one would expect the covering spectrum contains less information than the
length spectrum. This can also been seen in our example of a smooth one 
parameter family of nonisospectral tori with a common covering spectrum 
[Example~\ref{exdiamonds}].
Note that flat tori are isospectral iff they share the same length spectrum and
and are determined up to isometry by their marked length spectrum
\cite{Go3}. On the other hand length spectrum alone doesn't determine the
covering spectrum. We have examples of compact Riemannian 
manifolds
with a common length spectrum but a distinct covering spectrum 
[Example~\ref{samelap}]. 

In Section~\ref{sectmult} 
we define multiplicity [Defn~\ref{multi}] for the covering
spectrum, and find a bound on $\#_m(CovSpec(X)\cap [a,b])$ for $a>0$, where
$\#_m$ is the cardinality of
the set counting multiplicity [Lemma~\ref{tilNi(a,b)}].  We also define a
special set of generators of the revised fundamental group [Defn~\ref{multi}]
that we call the short basis. Roughly,
these are the elements of the revised fundamental group represented by loops
wrapped one time around a single hole in the space.  We prove this set
generates the revised fundamental
group in Proposition~\ref{mcpg} and show that the number of elements in this 
set
is $\#_m(CovSpec(X))$.

In Section~\ref{sectgh1} 
we focus on the relationship between the covering spectrum
and Gromov-Hausdorff convergence.  We begin by studying the Gromov-Hausdorff
convergence
of the delta covers, proving that if $X_i$ converge to $X$ in the GH sense, 
then
a subsequence of the delata covers $\tilde{X}_i^\delta$ converges as well
[Theorem~\ref{delta-precom}]. This involves
reworking Gromov's precompactness theorem and carefully
controlling the 
group of deck transforms of a delta cover.  In Example~\ref{subseq} we show it
is necessary
to use a subsequence and in Example~\ref{dimtoinf} we show that universal 
covers
need not have converging subsequences.  An immediate application is that if we
have a GH compact
class of compact length spaces with universal covers, then for $b>a>0$,
$\#_m(CovSpec(X)\cap [a,b])$ 
is uniformly bounded on this class [Cor~\ref{improved}]. 
One can also use the precompactness of the $\delta$-covers to show that the
{\em revised} fundamental groups of
such a compact class
with an additional uniform lower bound on the first systole, have finitely 
many
isomorphism classes extending Theorem 5 in \cite{ShW}.  

In Section~\ref{sectgh2},
we prove that if compact length spaces $X_i$ converge to a compact length
space $Y$ in the GH sense, then the covering spectra converge
[Thm~\ref{covspecconv}].  In particular,
\be
\lim_{i\to\infty} d_H(CovSpec(X_i)\cup \{0\}, CovSpec(Y)\cup \{0\}) \to 0
\ee
where $d_H$ is the Hausdorff distance between subsets of the real line
[Cor~\ref{corcovspconv}].  Note that
it is easy to see that when $1/j \times 1$ tori converge to a circle, there 
are
elements
of the covering spectrum which converge to $0$. Note also that if $M_i$ are
compact Riemannian manifolds
with $Ricci(M_i) \ge -(n-1)H$ and $diam(M_i)\le D$, such that $M_i$ converge 
to
$Y$, then $\#(CovSpec(M_i)) \ge \#(CovSpec(Y))$ for $i$ sufficiently large 
(not
counting multiplicity) [Cor~\ref{Ricci-number}].

We also prove that connected classes of compact length spaces with a common
discrete length spectrum, have a common covering spectrum
[Theorem~\ref{connlencov}].  In particular, a one parameter family of compact
Riemannian manifolds with a common length spectrum must have a common covering
spectrum [Cor~\ref{corconnlencov}].

In Section~\ref{sectgap} 
we study the gap phenomenon of the covering spectrum in certain
classes of compact
length spaces with universal covers [Prop~\ref{gap} and Prop~\ref{clumping}].
We apply these results
and \cite[Theorem 1.1]{SoWei}
to describe the gap and clumping properties of the covering spectra of
Riemannian manifolds
with $Ricci(M_i) \ge -(n-1)H$ and $diam(M_i)\le D$ and their limit spaces 
[Cor~\ref{corogap} and Cor~\ref{coroclump}].

In Section~\ref{sectlap} 
we relate the covering spectrum with the Laplace spectrum of a
manifold.  We first
easily show that if we have a class of negatively curved compact Riemannian
manifolds with a common Laplace spectrum, then there are only finitely
many possible covering spectra in this class [Prop~\ref{isofinite}].   
We conjecture that
this is true without the negative sectional curvature condition but with a
uniform upper bound on diameter [Conj~\ref{conjlap}].  
In Example~\ref{isolapnotcov}
we give a pair of
Heisenberg manifolds which are Laplace isospectral and yet have distinct
covering spectra.
This example heavily uses the work of Carolyn Gordon in \cite{Go1}, but it
should be noted that her
famous pairs of Laplace isospectral Heisenberg manifolds with distinct length
spectra in fact have
the same covering spectrum. 
We close by demonstrating that special pairs of Sunada isospectral manifolds,
the ones he attributes
to Komatsu \cite[Ex 3]{Su}, always share the same covering spectrum and, in
fact,
have only one element in
that covering spectrum [Prop~\ref{XKomatsu}].

The authors would like to thank Carolyn Gordon for
conversations and references regarding  the
length and Laplace spectra.  We'd also like to thank Jeff Cheeger for 
suggesting
that we investigate 
the possible existence of a gap phenomenon related to Theorem 1.1 of
\cite{SoWei}.

 \sect{Background}  \label{background}
 
 First we recall some basic definitions.
 
 \begin{defn} \label{compactlength} 
{\em
A {\em complete length space} is a 
complete metric space such that every pair of points in the space is joined by
 a length minimizing rectifiable curve.  The distance between
 the points is the length of that curve.  A {\em compact length
 space} is a compact complete length space.  (c.f. \cite{BBI}).
}
 \end{defn}
 
 \begin{defn} {\em
 We say $\bar{X}$ is a  {\em covering space} of $X$ if there
 is a continuous map $\pi: \bar{X} \to X$ such that $\forall x
 \in X$ there is an open neighborhood $U$ such that
 $\pi^{-1}(U)$ is a disjoint union of open subsets of
 $\bar{X}$ each of which is mapped homeomorphically onto $U$
 by $\pi$ (we say $U$ is evenly covered by $\pi$).   }
 \end{defn}
 
 Let $\cal U$ be any open covering of $Y$. For any $p \in Y$, by
 \cite[Page 81]{Sp}, there is a covering space, $\tilde{Y}_{\cal U}$,
 of $Y$ with covering group  $\pi_1(Y,{\cal U}, p)$, where
 $\pi_1(Y,{\cal U}, p)$ is a normal subgroup of $\pi_1(Y, p)$, generated
 by homotopy classes of closed paths having a representative of the form
 $\alpha^{-1} \circ \beta \circ \alpha$, where $\beta$ is a closed path
 lying in some element of $\cal U$ and $\alpha$ is a path from $p$ to
 $\beta(0)$.
 
 Now let's recall the $\delta$-covers we introduced in \cite{SoWei}.
 
\begin{defn} \label{defdel}
{\em
 Given $\delta > 0$, the  $\delta$-cover, denoted $\tilde{Y}^\delta$,
 of a length space $Y$, is defined to be $\tilde{Y}_{{\cal U}_{\delta}}$
 where ${\cal U}_\delta$ is the open covering of $Y$ consisting of
 all balls of radius $\delta$.
 
 The covering group will be denoted $\pi_1(Y,\delta, p)\subset \pi_1(Y,p)$
 and the group of deck transformations of $\tilde{Y}^\delta$ will
 be denoted $G(Y,\delta)=\pi_1(Y, p)/\pi_1(Y,\delta, p)$.
} 
\end{defn}
 
 It is easy to see that a delta cover is a regular or Galois cover.
 That is, the lift of any closed loop in Y is either always closed
 or always open in its delta cover.
 
 We now state some very simple lemmas:
 
 \begin{lem} \label{liftballs}
 If $\pi:\bar{Y} \to Y$ is a covering map between complete length
 spaces and $\forall y\in Y$, $\pi^{-1}(B_y(r))$ is a disjoint
 collection of balls of radius $r$ in $\bar{Y}$, then
 $\tilde{Y}^r$ covers $\bar{Y}$.
 \end{lem}
 
 \Pf
 Recall that in \cite [Ch2, Sec 5, Lem 11] {Sp}, Spanier shows that
 if $\pi: \bar{Y} \ra Y$ is a covering projection and $\cal U$
 is an open covering of $Y$ such that each of its open sets is
 evenly covered by $\pi$, then $\tilde{Y}_{\cal U}$ covers
 $\bar Y$.  Here $\cal{U}$ is the collection of balls of radius
 $r$, so we need only show that these balls are evenly covered by
 $\pi$.  
 
 Let $B_{\bar{y}}(r) \subset \pi^{-1}(B_y(r))$.  We need only
 show $\pi: B_{\bar{y}}(r) \ra B_y(r)$ is a homeomorphism.
 In fact, by the hypothesis, it is a covering map.  Thus if it
 is not 1:1, there are two preimages of $y$: $\bar{y}_1$ and $\bar{y}_2$.
 Note that $B_{\bar{y}_i}(r)$ is a connected subset of
 $\pi^{-1}(B_y(r))$, so it is a subset of $B_{\bar{y}}(r)$
 in which case $\bar{y}=\bar{y}_i$ and $\pi$ is 1:1.
 \qed
 
 \begin{ex}
 Suppose $Y$ is a flat $3\times2$ torus: $S^1_3\times S^1_2$, then it has the
 following delta covers:
 \begin{eqnarray*}
 \tilde{Y}^\delta &=& Y \textrm{ for } \delta > 3/2, \\
 \tilde{Y}^\delta &=& S^1_3 \times \mathbb R \textrm{ for } \delta \in (1, 
 3/2],
 \\
 \tilde{Y}^\delta &=& \mathbb R^2  \textrm{ for } \delta \in (0,1].
 \end{eqnarray*}
 \end{ex}
 
 \begin{lem}\label{mono} 
 The $\delta$ covers of complete length spaces are monotone
 in the sense that if $r<t$ then $\tilde{X}^r$ covers $\tilde{X}^t$.
 In fact $\tilde{X}^r$ is the $r$-cover of the complete length space
 $\tilde{X}^t$.
 \end{lem}
 
 \Pf
 Let $Y_t=\tilde{X}^t$.  We need only show $\tilde{Y_t}^r=\tilde{X}^r$
 for $r<t$.
 
 Applying Lemma~\ref{liftballs} to balls of radius $r$ in $X$.
 These must lift to unions of disjoint balls of radius $r$ in $\tilde{Y_t}^r$.
 Thus $\tilde{X}^r$ covers $\tilde{Y_t}^r$.
 
 Recall from \cite[Ch2, Sec 5, 8]{Sp} that if $\cal V$ is an open covering of 
 $X$ that refines $\cal U$, then
 $\pi_1(X,{\cal V}, p) \subset \pi_1(X,{\cal U}, p)$, or $\tilde{X}_{\cal V}$
 covers $\tilde{X}_{\cal U}$.  
 
 Thus, clearly, $\tilde{X}^r$ covers $Y_t=\tilde{X}^t$.  
 
 Now we apply Lemma~\ref{liftballs} to balls of radius $r$ in $Y_t$.
 These must lift to unions of disjoint balls of radius $r$ in $\tilde{X}^r$,
 as can be seen by first projecting them down to $X$.  Thus
 $\tilde{Y_t}^r$ covers $\tilde{X}^r$ and we are done.
 \qed
 
 In the following lemma we restrict ourselves to compact length spaces.
 
 \begin{lem} \label{lower}
 The $\delta$-covers of a compact length space $X$ are lower semi-continuous.
 In fact, for any $\delta>0$, there exists $\epsilon\in (0, \delta)$ such that
    $\tilde{X}^\epsilon=\tilde{X}^\delta$.
 \end{lem}
 
 \Pf If not, there is a sequence of $\delta_i \ra \delta$ increasingly, such 
 that
 $\tilde{X}^{\delta_i}\not = \tilde{X}^\delta$ for each $i$. Namely there 
exist 
 a
 sequence of closed curves $\gamma_i$ in $X$ with length $l(\gamma_i) \leq 
 2\diam
 (X) + 2\delta_i$, which lifts to an open curve in $\tilde{X}^{\delta_i}$ but 
 an
 closed curve in $ \tilde{X}^\delta$. Parametrize each curve
 by the unit interval $[0,1]$ with constant speed. Since $X$ is
 compact, by Arzela-Ascoli
 theorem, there is a subsequence of $\gamma_i$ which
 converges to some closed curve $\gamma: [0,1] \ra X$ uniformly. So 
 $d(\gamma_i(t),\gamma(t)) < \delta/2$ for all $i$ large and $t \in [0,1]$.
 Hence $\gamma_i, \gamma$ lift the same to the covering spaces
 $\tilde{X}^{\delta_i},\ \tilde{X}^{\delta}$ for $i$ large.  That is, $\gamma$
 lifts to an open curve in $\tilde{X}^{\delta_i}$ for
 all $i$ large and a closed curved in $ \tilde{X}^\delta$. From
 Definition~\ref{defdel}, $\gamma$ lies in some
 finite union of open $\delta$-balls in $X$, so it must also lie in some 
union 
 of open
 $\delta'$-balls for some $\delta' < \delta$, which contradicts to that 
 $\gamma$ 
 lifts to an open curve in $\tilde{X}^{\delta_i}$ for all $i$ large.
 \qed
 
 \begin{ex}
 The Hawaii ring (c.f. \cite{Sp}) is a compact length space which consists
 of an infinite set of rings of radii $r_i$ decreasing to $0$, all
 joined at a common point.  This space has an infinite sequence of
 distinct $\delta$ covers as $\delta$ converges to $0$.
 
 A modified complete noncompact Hawaii ring can be created by
 taking an infinite set of rings of radius $r_i$ increasing to
 $r_0=1$.  This space also has an infinite sequence of distinct $\delta$
 covers as $\delta$ approaches $r_0$.  
 
 In both of these spaces, the $\tilde{X}^{\pi r_i}$ are all distinct
 covers.
 \end{ex}
 
 This example demonstrates
 that the compactness hypothesis in Lemma~\ref{lower} is necessary.

 %**************************************************************
 
 \sect{The covering spectrum} \label{delcovsect}
 
 We now define the covering spectrum by singling out the 
 deltas where the delta covering spaces change.
 
 \begin{defn}  \label{coveringspectrum}
 Given a complete length space $X$, the covering spectrum of
 $X$, denoted  CovSpec$(X)$ is the set of all $\delta > 0$ such
 that
 \be
 \tilde{X}^\delta \neq \tilde{X}^{\delta'}
 \ee
 for all $\delta'> \delta$.
  \end{defn}
 
 Since the $\delta$-covers are monotone, this is equivalent
 to say for any $\epsilon > 0$, there exists $\delta'$ with
 $0<\delta'-\delta<\epsilon$ such that $$ \tilde{X}^\delta
 \neq \tilde{X}^{\delta'}. $$
 
 In general, for a compact length space $X$, the CovSpec$(X)$
 lies in $(0, \diam (X)/2]$.
 
 In our above examples the
 covering spectrum of the flat $3$ by $2$ torus is $\{1,3/2\}$,
 and the traditional Hawaii ring with infinite circles of radii $r_i$ is
 $\{\pi r_i : i \in \mathbb N \}$.
 
 We have the following property of the covering spectrum.
 
 \begin{prop} \label{discrete}
 For a compact length space, $X$, its CovSpec$(X)$ is discrete
 and  $Cl(CovSpec(X) \subset \mathbb R)\subset CovSpec(X)\cup \{0\}$. 
 \end{prop}

 \Pf Since zero is not in CovSpec$(X)$, if CovSpec$(X)$ is not discrete,
  we can assume it has an accumulation at some $\delta > 0$. In fact we can 
 assume there is a strictly decreasing sequence of $\delta_i \in CovSpec (X)$
 converges to the $\delta > 0$ since $\delta$-covers are lower semi-continuous
 (Lemma~\ref{lower}).
 Let $\gamma_i$ be a loop at $p \in X$ such that $\gamma_i$
 lifts trivially to $\tilde{X}^{\delta_i}$ but
 nontrivially to $\tilde{X}^{\delta_{i+1}}$ with length $l(\gamma_i) \leq 
 2\diam
 (X) + 2\delta_i$. 
Parametrize each curve in the sequence 
by the unit interval $[0,1]$ with constant speed.
Since $X$ is compact, by the Arzela-Ascoli
 theorem, there is an uniformly converging subsequence, which we will
still call $\gamma_i$.
  So there is an $i$ sufficiently large  that $d(\gamma_i(t),
 \gamma_{i+1}(t)) < \delta$ for all $t \in [0,1]$. Since $\delta_i > \delta$
 the covering maps are isometric on $\delta$ balls for all $i$, so
 $\gamma_i, \gamma_{i+1}$ lift the same to the covering spaces
 $\tilde{X}^{\delta_i},\ \tilde{X}^{\delta_{i+1}}$, contradicting that
 $\gamma_{i+1}$ lifts trivially
 to $\tilde{X}^{\delta_{i+1}}$ and $\gamma_i$ lifts
 nontrivially to $\tilde{X}^{\delta_{i+1}}$. Therefore CovSpec$(X)$ is 
 discrete.
 \qed

 The example of the Hawaii ring shows that $0$ could be in
 the closure of the covering spectrum of a compact length space. 
 Proposition~\ref{discrete} is not true for a noncompact 
complete length space as 
 another revised Hawaii ring,
 the union of the sequence of circles with a common point
 and radius $r_i$ decreasing to $1$, shows.
 
 We now turn to a discussion of the existence of universal covers.
 The original compact Hawaii ring with $r_i \to 0$
 is a classic example of a compact length space with no universal cover.
 Recall the definition of a universal cover.
 
 \begin{defn}\label{univcov} \cite[pp
 62,83]{Sp} {\em We say $\tilde{X}$ is a
 {\em universal cover} of $X$ if $\tilde{X}$ is a cover of $X$
 such that for any other cover $\bar{X}$ of $X$, there is a
 commutative triangle formed by a continuous map $f:\tilde{X}
 \to \bar{X}$ and the two covering projections.} 
 \end{defn}
 
 In \cite[Prop. 3.2]{SoWei} we proved that if a compact length
 space $Y$ has a universal cover $\tilde{Y}$ then $\tilde{Y}$ is a delta 
cover. 
 In fact, $Y$ has a universal cover iff the delta covers
 stabilize: there exists a $\delta_0 > 0$ such that
 $\tilde{Y}^\delta=\tilde{Y}^{\delta_0}$ for all
 $\delta<\delta_0$ \cite[Thm 3.7]{SoWei}.
 Clearly the delta covers of the Hawaii ring do not stabilize.

  \begin{theo} \label{univ} For a compact length space $X$,
 its universal cover $\tilde{X}$ exists iff it's covering
spectrum,  CovSpec$(X)$,
 is finite.
  \end{theo}
 
 \Pf If the CovSpec$(X)$ is finite, then $\delta_0 = \min \{
 CovSpec (X) \}$ is positive. So the $\delta$-covers
 stabilize and
by \cite[Thm 3.7]{SoWei} the universal
 cover of $X$ exists.
 
 If the universal cover ${\tilde X}$ of $X$ exists, then
 CovSpec$(X)$ lies in $[\delta_X, \diam (X)/2]$ for some
 $\delta_X > 0$. By Proposition~\ref{discrete} 
the CovSpec$(X)\cap [\delta_X, \diam (X)/2]$ is closed
and discrete.   Therefore CovSpec$(X)$ is finite.
 \qed
 
 Although the covering spectrum is defined using layers
 of covering spaces, $\#\{CovSpec(X)\}$ does not count the
 total number of covering spaces of $X$.  The $\delta$
 covers are a very small selection of covering spaces.  Clearly
 the $3\times2$ torus has many covering spaces that are tori 
 and cylinders which are not delta covers (i.e.
 $S_{3i}^1 \times S_{2j}^1$, $\mathbb R \times S_{2j}^1$ 
and $S_{3j}^1\times \mathbb R$).  
 Furthermore, the lens spaces,
 $S^3$ mod $\mathbb Z^k$ with the standard metric only have one 
$\delta$-cover,
 $S^3$, although they often have many covering spaces.
 
 The covering spectrum can intuitively be though of as capturing the
 size of holes in the length space.  For the $2\times3$ torus, it captures
 information about both of the holes in the torus: both of the
 generators of the fundamental group.  The fact that $\#\{CovSpec(X)\}=2$ 
 in this example is strongly related to the fact that there are two
 generators of the fundamental group.
 
 On the other hand, the covering spectrum of a $1\times 1$ torus has only
 one element because both holes in this torus have the same size.
 Later on, we will define multiplicity for the
 elements of the covering spectrum, which will better enable us
 to capture the fact that there are two ``holes'' in this torus as well.
 
 %**************************************************************************
 \sect{The covering spectrum and length spectrum}  \label{covandlength}
 
 In this section we restrict our attention to complete length spaces $X$ 
 which have a universal cover. 
 
 First recall that a geodesic in a length space is a curve which is locally a 
 distance minimizer in the following sense \cite{BBI}. 
 
 \begin{defn}  \label{geodesicdef}
{\em  A curve $\gamma: I \ra X$ is
 called a {\em geodesic} if for every $t \in I$ there exists an interval $J$
 containing a neighborhood of $t$ in $I$ such that $\gamma|_{J}$ is a shortest
 path.  A {\em closed geodesic} is a geodesic loop which is minimizing
in a neighborhood of its end point.}
 \end{defn}
 
 It is easy to use the definition of a covering space to show that a length
minimizing curve in a covering space projects to a geodesic,
 and that geodesics lift to geodesics.
 
 Then one can naturally extend the definition of length spectrum from
 manifolds to complete length spaces.
 
 \begin{defn}  \label{lsdef}
{\em 
The {\em length spectrum}, Length($X$), of a complete length space, $X$,
 is the set of lengths of
 closed geodesics.  It is counted with multiplicity where 
 the multiplicity refers
 to the number of distinct free homotopy classes that contain a closed 
 geodesic of that length.  }
 \end{defn}
 
We recall the definition of the revised fundamental group from \cite{SoWei}.

 \begin{defn}\label{revfg}  {\em
 The {\em revised fundamental group}, $\bar{\pi}_1(X)$, 
 of a complete length space, $X$, with
 a universal cover, $\tilde{X}$,
is the group of deck transforms of the universal cover.  
 Given an element, $g\in \bar{\pi}_1(X)$, and a base point $x\in X$
 a {\em representative loop} of $g$ based at $x$ is a curve, $c$,
such that $c(0)=x$ whose lift $\tilde{c}$ to the universal cover
 runs from a point $\tilde{c}(0)$ to $g\tilde{c}(0)$.
 }
 \end{defn}
 
For simplicity the reader may wish to assume $X$ has a simply connected 
 universal cover, or equivalently, that $X$ is  semi-locally simply connected.
 In that case, the fundamental group $\pi_1(X)$ of $X$ is isomorphic
 to $\bar{\pi}_1(X)$.  

 In general, however, the universal cover of a compact length space
 may not be simply connected.
 One example is the double suspension over the Hawaii
 Ring (c.f. \cite{Sp}) which is its own universal cover but has an
 infinite fundamental group because the infinite alternation of loops in the
Hawaii rings aren't contractible.  These loops are homotopic to  loops in an
arbitrarily small neighborhood but not to a single point. 

When the universal cover is not simply connected, the representive
loops of the identity element are the projections of arbitrary loops
in the universal cover, which are not necessarily contractible.  Thus the
equivalence class of representative loops corresponding to an
element $g\in \bar{\pi}_1(X)$ and a point $x\in X$
is not a homotopy equivalence class,
but rather a collection of homotopy equivalence classes.  

The following lemma is easy to prove using the fact that the universal
 cover is a $\delta$-cover and using the compactness of $X$.
 
 \begin{lemma}  \label{shortest}
 Given a compact length space $X$ with a universal cover $\tilde{X}$, 
 for all nontrivial $g \in \bar{\pi}_1(X)$, we have
 \be
 m(g):=\min_{\tilde{x}\in \tilde{X}} d_{\tilde{X}}(\tilde{x}, g\tilde{x})
\subset Length(X).
 \ee
 If $\gamma_g$ is the projection of a minimizing curve joining 
 a minimizing pair of points $\tilde{x}$ and $g\tilde{x}$, then $\gamma_g$ is a
closed 
geodesic 
 in $X$ of length $m(g)$ which is a shortest
 curve in its free homotopy class.%and R-homotopy class
 \end{lemma}
 
\Pf
There exists $\delta>0$ such that the universal cover is a $\delta$ cover.
Thus
\be
m(g)=\inf_{\tilde{x}\in \tilde{X}} d_{\tilde{X}}(\tilde{x}, g\tilde{x})\ge
2\delta>0.
\ee
Let $\tilde{x}_i \in \tilde{X}$ approach this infimum.  Since $X$ is compact,
a subsequence of $x_i=\pi(\tilde{x}_i)$ converges to some $x$ whose lift
$\tilde{x}$ then achieves this infimum.  Then $\gamma_g$ is the shortest
representative of $g$ for any base point, it has length $m(g)$
and it is the projection of a geodesic to a loop.  Extending the definition
of $\gamma_g$ periodically, we can see that it is a representative of $g$
based at $\gamma(t)$ as well. So it must be the projection of a length
minimizing curve between $\tilde{\gamma}_g(t)$ and $g\tilde{\gamma}_g(t)$ which
implies that it is a closed geodesic.
\qed

Thus we have the following useful map.

 \begin{defn}\label{minmarkmap}
{\em
 The {\em minimum marked length map} of a compact length space $X$
 with a universal cover is the function 
 $m:\bar{\pi}_1(X)\to LengthSpec(X) \cup \{0\}$  defined in 
 Lemma~\ref{shortest}.} 
 \end{defn}
 
\Rk: 
The minimum marked length map is closely related to the 
translative delta length
 $l(g,\delta)$ we defined in \cite[Definition 3.2]{SoWei}. Recall  
\be
l(g, \delta)=
\min_{q\in \tilde{X}^\delta}
d_{\tilde{Y}^\delta}(q, g(q)).
\ee
and $l(g,\delta) \ge 2\delta$ for all $g$ acts nontrivially on
$\tilde{X}^\delta$. Note that,
since covering maps are distance decreasing,
\be  \label{l-delta}
m(g) \ge l(g,\delta) \  \mbox{for all}\ g \in \bar{\pi}_1(X).
\ee
 
 \begin{lemma} \label{mfinite}
 When $X$ is compact, the set $Im(m)=\{m(g): g\in \bar{\pi}_1(X)\}$ is closed 
 and
 discrete.  Furthermore $m(g)=0$ iff $g=e$. 
 \end{lemma}
 
 \Pf
 First note that, by the Arzela-Ascoli Theorem, sequences of length 
 minimizing curves have subsequences which
 converge to length minimizing curves, so if $m_i \in Im(m)$ converge
 to $m_\infty$ then we have a subsequence of  $x_i$ in the fundamental 
 domain of $\tilde{X}$ converging to $x_\infty$, and $g_ix_i$ converging
 to some $y_\infty$ such that $d_{\tilde{X}}(x_\infty, y_\infty)=m_\infty$.
 Since the universal cover is a delta cover, $y_\infty=y_i$ for $i$ 
 sufficiently 
 large.  Thus $m_\infty=m_i$ for all $i$ large and $Im(m)$  is closed and
 discrete.
 
 We know $m(g)=0$ iff $g=e$ because only a trivial deck transform fixes a 
 point.
 \qed
 
 We have the following nice relation.
 
 \begin{theo} \label{cov-length-spectrum}
 When $X$ is a compact length space with a universal cover then
 \be
   2 CovSpec(X) \subset Im( m(\bar{\pi}_1(X) ) \subset LengthSpec(X) \cup 
 \{0\},
 \ee
 where $m$ is the minimum marked length map defined in
 Definition~\ref{minmarkmap}.
 \end{theo}
 
 This theorem follows from the following definition and lemma.
 
 \begin{defn}\label{deltapair}
{\em 
 If $X$ is a complete length space and $\delta>0$ then we say a 
 $\delta$-pair is a pair of points $\{x_1, x_2\}$ in $\tilde{X}^{\delta}$
 which are not equal but are projected to the same point in 
 $\tilde{X}^{\delta'}$
 for all $\delta'>\delta$.  }
 \end{defn}

 \begin{lemma} \label{lemdelpair}
 Fix a compact length space with a universal cover, $X$, and  $\delta\in
 CovSpec(X)$.
 Let $h_\delta=\inf d_{\tilde{X}^\delta}(x_1, x_2)$ over all $\delta$-pairs 
$x_1, x_2$.   Then this infimum is achieved, there is an element
 $g\in \bar{\pi}_1$ such that $m(g)=h_\delta$ and $h_\delta=2\delta$.
 \end{lemma}
 
 \Pf
 By compactness it is easy to show that there exists a $\delta$ pair $x_1, 
x_2$
 which achieves this infimum.  It is not necessarily a unique pair even up
 to deck transforms.
 
 First $h_\delta \geq 2\delta$, else a minimizing curve from
 $x_1$ to $x_2$ would have length $< 2\delta$ and its projection to $X$
 would fit in $B_{\pi(x_1)}(\delta)$, so it would be lifted as a loop
 to $\tilde{X}^\delta$ making $x_1=x_2$ by Definition~\ref{defdel}.
 
 Now we will show $h=h_\delta \le 2\delta$. By Proposition~\ref{univ}
		 the covering spectrum is finite, so
 there is $\epsilon > \delta$ such that for all $\delta' \in ( \delta,
 \epsilon)$,
 $\tilde{X}^{\delta'}=\tilde{X}^\epsilon$. Naturally $\tilde{X}^{\delta}$ is a
nontrivial cover of $\tilde{X}^\epsilon$.

 Note that $y_1$ and $y_2$ are a $\delta$ pair iff they are not equal but
 project to the same point in $\tilde{X}^\epsilon$.
 So we have for all $y_1\neq y_2 \in \tilde{X}^\delta$ such that
 $\pi_\epsilon(y_1)=\pi_\epsilon(y_2) \in \tilde{X}^\epsilon$,
 ${B}_{y_1}(h/2)$ and ${B}_{y_2}(h/2)$ are disjoint.
 So for all $z\in \tilde{X}^\epsilon, B_z(h/2)$ lifts to
 a disjoint union of balls in $\tilde{X}^\delta$.  Thus
 by Lemma~\ref{liftballs}, applied to $\bar{Y}=\tilde{X}^\delta$
 as a cover of $Y=\tilde{X}^\epsilon$, we get
 $\tilde{Y}^{h/2}$ covers $\bar{Y}=\tilde{X}^\delta$. If $h/2>\epsilon$, then
$\tilde{Y}^{h/2} = \tilde{X}^\epsilon$, which means $\tilde{X}^\delta =
\tilde{X}^\epsilon$.  This is a contradiction. So we have $h/2 \le \epsilon$,
then
$\tilde{Y}^{h/2} = \tilde{X}^{h/2}$, so $\tilde{X}^{h/2}$ covers
$\tilde{X}^\delta$. Therefore $h\le 2\delta$.
 
 So $h_\delta=2\delta$.
 
 Now let $C$ be a minimal geodesic connecting $x_1$ and $x_2$ and  $g$ an
element in $\bar{\pi}_1(X)$ which is represented by the projection of $C$. 
Then $m(g) \le 2\delta$. But $g$ acts
nontrivially on $\tilde{X}^\delta$, so $l(g,\delta) \ge 2\delta$. Therefore
$m(g) \ge 2\delta$ by (\ref{l-delta}). Thus $m(g) = 2\delta$.
 \qed

Another standard length spectrum defined on manifolds is the
minimal length spectrum. 
 
 \begin{defn} \label{mlsdef}{\em
The {\em minimal length spectrum} is the set of lengths of
 closed geodesics which are the shortest in their free homotopy class.}
  \end{defn}

If a compact length space $X$ is a semilocally simply connected, 
or equivalently
has a simply connected universal cover, then the above definition
makes sense and each homotopy class contains a curve of minimum
length.  In fact, the minimal length spectrum agrees with
$im(m)\setminus \{0\}$ as can be seen in the following lemma
combined with Lemma~\ref{shortest}.

\begin{lemma}\label{minsurjl}
For a compact length space with a simply connected universal cover
the minimum marked length map $m$ maps surjectively onto the 
minimal length spectrum $\cup \{0\}$.  
\end{lemma}

\Pf
Given any $L$ in the minimal length spectrum,
there is a free homotopy class of loops whose minimum length
is $L$.  Let $c_1$ be the shortest such loop. It defines a deck 
transform $g$ and $m(g)\le L(c)=L$.  
Suppose $m(g)<L$, then there exists $\tilde{x}\in \tilde{X}$
such that $d(g\tilde{x}, \tilde{x})<L$.  Join this pair of
points by a length minimizing curve $\tilde{c}_2$.  

% adjust appropriately
If
the universal cover is simply connected then the projection
$c_2$ is a loop freely homotopic to $c_1$ and we have a contradiction.
\qed

Theorem~\ref{cov-length-spectrum},
Lemma~\ref{mfinite} and Lemma~\ref{minsurjl}
combine to give us the following theorem:

 \begin{theo} \label{2cov-length-spectrum}
 When $X$ is a compact length space with a simply
connected universal cover then the minimum
length spectrum is closed and discrete and
 \be
   2 CovSpec(X) \subset MinLengthSpec(X).
 \ee
 \end{theo}

%***************************************
\sect{The marked length spectrum}
\label{markcovandlen}

 A stronger concept than the length spectrum of a manifold
 is the marked length spectrum which includes information
 about the fundamental group itself.  Here we will study
arbitrary compact length spaces with universal covers.
The natural extension
 of the definition of marked length spectrum to such spaces
 involves the revised fundamental group $\bar{\pi}_1(X)$ 
instead of the fundamental group
[Defn~\ref{revfg}].  For simplicity, the reader may wish to
assume the universal cover is simply connected in which
case the revised fundemantal group is just the fundamental
group of the space.
 
 \begin{defn}  \label{markdef}{\em
 Given a complete length space $X$, the {\em marked length spectrum} 
 of $X$ is a function $MLS$ that
 associates to each element $g$ in $\bar{\pi}_1(X)$ the set of lengths, 
$MLS(g)$, of the
 closed geodesics freely homotopic to a representative loop of $g$.
 Clearly, this map only depends on the conjugacy class of $g$.
 
 Two spaces $X_1$ and $X_2$ are said to have the same marked length
 spectrum iff there is an isomorphism between their
 revised fundamental groups which commutes with their marked
length maps $MLS_1$ and $MLS_2$.}
\end{defn}
 
Recall the definition of the minimum marked length map,
$m:\bar{\pi}_1(X)\to (0,\infty)$, in
Definition~\ref{minmarkmap} and Lemma~\ref{shortest}.  
Since $MLS(g)$ includes the lengths of all geodesics 
representing $g$, we have $m(g)=min(MLS(g))$.

\begin{defn}\label{minmarklspec}{\em
 We say two spaces with universal covers $X_1$ and $X_2$ have the same 
 {\em minimum marked length spectrum} 
 iff there is an isomorphism between their
 revised fundamental groups which commutes with their minimum marked
length maps $m_1$ and $m_2$.}
 \end{defn}

  We can also mark the covering spectrum of a compact length space
 with a universal cover using the following simple map:
 
 \begin{defn}\label{covspecmap}
 {\em Given a complete length space $X$ with a universal cover,
 we define the {\em covering spectrum map},
  $f:\bar{\pi}_1(X) \to CovSpec(X) \cup \{0\} $ as follows.
 
 Given $g \in \bar{\pi}_1(X)$, let $f(g)$ be the unique $\delta$ in CovSpec(X)
 such that a loop $\gamma_g$ representing $g$ in $X$ lifts to a curve in 
 $\tilde{X}^\delta$
 that is not a loop, but lifts to a loop in $\tilde{X}^{\delta'}$ for all
 $\delta'>\delta$.  In fact all loops freely r-homotopic to this one
 will then also share this property.
 
Equivalently $f(g)$ is the largest $\delta$ such that the projections
 of $gx$ and $x$ from $\tilde{X}$ to $\tilde{X}^\delta$ are distinct points.
}
 \end{defn}
 
 Note $f(g)=0$ iff $g=e$.

 \begin{lem}  \label{surj}
 Given a compact length space $X$ with a universal cover,
 the covering spectrum map $f:\bar{\pi}_1(X) \to CovSpec(X) \cup \{0\} $ 
 is surjective.
 \end{lem}
 
 \Pf
 If $\delta\in CovSpec(X)$, then $\tilde{X}^\delta\neq \tilde{X}^{\delta'}$
 for all $\delta'>\delta$.  Since $X$ is compact, the covering spectrum
 is discrete away from 0 [Theorem~\ref{discrete}], so there exists
$\epsilon>\delta$ such that
 $\tilde{X}^\epsilon=\tilde{X}^{\delta'}$ 
 for all $\delta'\in (\delta, \epsilon]$.
 Let $x_1$ and $x_2$ be a pair of distinct points in $\tilde{X}^\delta$
 which are mapped to the same point in $\tilde{X}^\epsilon$.  Let $C$
 be a curve joining $x_1$ to $x_2$.  Then $C$  projects to a loop in
 $X$, which lifts as a loop to $\tilde{X}^{\delta'}$ for all $\delta'>\delta$
and lifts to a curve that is not a loop in $\tilde{X}^\delta$.  
Let $g\in \bar{\pi}_1(X)$ which is represented by the projection of $C$, then
$f(g)=\delta$.  
 \qed
 
 Note that the compactness in this lemma is necessary as the following example 
 shows. Let $X$ be a revised Hawaii ring with circles of radius $1+1/n$ all
 attached at one point, then $X$ has a universal cover.
 Furthermore, $\pi$ is in CovSpec($X$), but it doesn't lie
 in the image of $f$ since the circle of radius $1$ is not in $X$.

 \begin{lem}  \label{subgroup}
 When $X$ is a complete length space with a universal cover
 and $f$ is the covering
 spectrum map, then 
 $f^{-1}([0,\delta])$ is a subgroup of $\bar{\pi}_1(X)$.
 \end{lem}
 
 \Pf
 If $g_1, g_2 \in f^{-1}([0,\delta])$, suppose $f(g_i)=\delta_i$, then there 
 are
 loops $\gamma_i$ which
 lifts as closed loop to $\tilde{M}^{\delta' }$ for $\delta'>\delta_i$ and as 
 an
 open curve to $\tilde{M}^{\delta_i}$.  So for any $\delta'>\max
 \{\delta_1,\delta_2\}$ both curves $\gamma_1, \gamma_2$ lift as closed loops 
 to
 $\tilde{M}^{\delta' }$. Now the element $g_1g_2$ can be represented
 by the loop $\gamma_1$ following $\gamma_2$, so the lift of the combination 
is
 closed in $\tilde{M}^{\delta' }$.  Thus 
 \be \label{max}
 f(g_1g_2)\leq \max \{\delta_1,\delta_2\}=\max\{f(g_1),f(g_2)\}.
 \ee
 \qed

 A nonpositively curved metric on a surface of genus $\ge 2$ with the set 
where
the curvature is $0$ has empty interior is determined
 up to isometry by its marked length spectrum \cite{Ot, Cr, Fa}. The same is
true for flat tori \cite{Go3}. 
 The following example demonstrates that even on flat tori, the covering
 spectrum does not determine the isometry class.  In fact it includes a smooth
 family of flat tori with a common covering spectra.
 
 \begin{ex} \label{exdiamonds}
 Here we examine a set of flat 2 dimensional tori, $T^2_\theta$, defined
 as rhombi with side length 1 and a variable angle $0<\theta \le 
\frac{\pi}{2}$
 between the sides.  Opposite sides are identified in the
 usual way and the universal cover of any of these examples
 is the Euclidean plane.  Note that, in this case, the
 marked length spectrum has only one length per element
 of the abelian fundamental group, so we can denote it
 as $m(g)$.
 
 If we locate the fundamental domain with corners at $(0,0)$, $(1,0)$,
$(cos(\theta),sin(\theta))$
 and $(1+cos(\theta), sin(\theta))$ then the group of deck transforms is
generated by $g_1: (x,y) \mapsto (x+1,y)$
 and $g_2:(x,y) \mapsto (x+cos(\theta), y+sin(\theta))$.
 
 Now for $\theta \in [\pi/3, \pi/2]$ it is easy to see that
 $m(g_1^a g_2^b)=\sqrt{a^2+b^2+2abcos(\theta)}$ so the
 length spectrum is 
\be
\{\sqrt{a^2+b^2+2abcos(\theta)}: a,b \in \mathbb Z \setminus \{0\} \}.
\ee
 Furthermore  $f(g_1^ag_2^b)=1/2$ unless $a=b=0$, so the covering spectrum
 is just $\{1/2\}$.  This provides us with a one parameter
 family of flat tori with a common covering spectrum.
 
 Now for $\theta \in (0, \pi/3)$, we get the same formula
 $m(g_1^a g_2^b)=\sqrt{a^2+b^2+2abcos(\theta)}$ so the
 length spectrum is 
\be
\{\sqrt{a^2+b^2+2abcos(\theta)}: a,b \in \mathbb Z\setminus\{0\} \}.
\ee
 However, now 
  \[
 f(g_1^ag_2^b) = \left\{  \begin{array}{ll} (1/2)\sqrt{1+1-2cos(\theta)} < 1/2 
 &
 \textrm{ if } a=-b \\
  1/2 & \textrm{ otherwise}
 \end{array} \right. .
 \]
 So the covering spectrum has 2 distinct elements.
 
 Since these two families together form a single one parameter family,
 we have also shown that the number of elements of the covering
 spectrum may change.  Although it is nice to see that here the
 covering spectra do vary continuously in Hausdorff sense.  
 \end{ex}

 In fact the covering spectrum is determined by the minimum marked
 length spectrum [Defn~\ref{minmarklspec}].
 
 \begin{theo}  \label{mlscs}
 Let $X_1$ and $X_2$ be compact length spaces
 with universal covers.
 If they have the same minimum marked length
 spectrum then they have the same covering spectrum.
 
 In fact if $CovSpec(X)=\{\delta_1< \delta_2 < \cdots <\delta_k\}$ then
 there exists a special sequence of subgroups
 $\{e\}=G_0 \subset G_1 \subset G_2 \subset \cdots G_k= \bar{\pi}_1(X)$ 
 such that each $G_i$ is generated by 
 \be \label{mlscseq}
 S_i=\{h\in \bar{\pi}_1(X): m(h)=2\delta_i\}
 \ee 
 combined with the elements of $G_{i-1}$.  Furthermore
 $f(g)=\delta_i$ implies $g\in G_i$ while $g\in G_i\setminus G_{i-1}$
 has $f(g)=\delta_i$.
 \end{theo}
 
 First we state a simple lemma which we will need.
 
 \begin{lem} \label{balltoL}
   Suppose $C:[0,L] \to B_q(\delta)\subset X$ where $X$ is a complete length
 space,
 then $C$ is freely homotopic
 to a product of curves of length $<2\delta$ based at $q$.
 \end{lem}
 
 \Pf
 We assume $C$ is parametrized by arclength.  Since its image is closed,
 it is in fact contained in $B_q(\delta-\epsilon)$ for some $\epsilon>0$.
 Partition $[0,L]$ into pieces of length $<\epsilon$: 
$t_1=0<t_2<t_3<...<t_k=L$
 Let $\sigma_i$ run minimally from $q$ to $C(t_i)$ so it has length
 $<\delta-\epsilon$.  Set $\sigma_k=\sigma_1$.
 So $C_i$ starting at $q$ running along $\sigma_i$
 to $C(t_i)$ running along $C$ to $C(t_{i+1})$ and running backwards
 along $\sigma_{i+1}$ to $q$ is a closed curve of length $< 2\delta$
 
 The product of these $C_i$ is a curve which is freely homotopic to $C$
 (where the homotopy runs along $\sigma_1=\sigma_k$.).
 \qed
 
 \begin{coro}  \label{L2del}
 If $C:[0,L]\to B_q(\delta) \subset X$ parametrized by arclength is the 
 shortest
 noncontractible curve in $X$, then $L<2\delta$.
 \end{coro}
 
 \noindent
 {\bf Proof of Theorem~\ref{mlscs}:}
  We will derive the marked covering map $f: 
 \bar{\pi}_1(X) \ra CovSpec(X) \cup \{0\}$ from the 
 marked shortest length spectrum $m$.
 
 We first claim that 
  \be  \label{flef1}
  f(g) \le (1/2)m(g) \ \  \mbox{ for all} \ g \in \bar{\pi}_1(X).
  \ee
 
 If $(1/2)m(g)<f(g)$ then a representative of $g$ is freely r-homotopic to a 
 curve of length $< 2f(g)$.  Such a curve must be contained in a ball of 
 radius $f(g)$, so it would lift as a closed curve to the $f(g)$ cover
 of $X$, but it cannot by the defn of $f(g)$ [Defn~\ref{covspecmap}].
 
 Let $G_0=\{e\}$ and $\delta_0=0$. We will construct the covering spectrum and 
 the map $f$ by induction.
 
  Suppose 
  
  a) we've defined distinct subgroups
 $G_0\subset G_1 \subset  \cdots G_k\subset \bar{\pi}_1(X)$ and defined
 $\delta_0 <\delta_1 < \cdots <\delta_k$ such that each
  \be
  \delta_i = min \{ m(g): g\in \pi_1 \setminus G_{i-1} \} \subset \{0\} \cup
 CovSpec(X).
  \ee
 (Here $G_{-1}$ is the empty set.)
   
   b) each $G_i$ is generated by all the
  elements $h\in \pi_1$ such that $m(h)=\delta_i$ combined with
  the elements of $G_{i-1}$.  
  
  c) each $G_i$
  contains every element $g\in \bar{\pi}_1$ such that $f(g)=\delta_i$
 and
 \be  \label{deltaspec}
 \{\delta_0,\delta_1, \cdots, \delta_k\}
 =CovSpec(X)\cap [0,\delta_k] \cup \{0\}.
 \ee
 
 Note that for all $g\in G_{i}\setminus G_{i-1}$ we have $f(g)=\delta_i$
 as a consequence of the above hypothesis b) c) and  (\ref{flef1}) 
(\ref{max}).
 
 The hypothesis a) b) c) are trivially true for $k=0$.
 
 Now we prove the induction step:
  
  If $G_k \neq \pi_1$, let
 \be  \label{eqdefdel}
 \delta_{k+1} = min \{(1/2) m(g): g\in \bar{\pi}_1 \setminus G_{k} \}.
 \ee
 Since $\{m(g)\}$ is closed and discrete [Lemma~\ref{mfinite}], 
 the minimum exists and is achieved.
 
 We will first show $\delta_{k+1}\in CovSpec(X)$.
 
 We have $h\in \pi_1 \setminus G_k$ such that $m(h)=2\delta_{k+1}$.
 We need to show $f(h)=\delta_{k+1}$.  By  (\ref{flef1})  $f(h)\le 
 \delta_{k+1}$.
  Assume on the contrary that $f(h)<\delta_{k+1}$.   
 Then for any $f(h) <\delta'<\delta_{k+1}$
 $h$ lifts trivially to $\tilde{X}^{\delta'}$.

 So $h$  is a product of elements $g_1g_2...g_k$ where
 each $g_k$ has a representative curve based at $p$ of the form
 $\alpha_i^{-1}\beta_i \alpha_i$ where $\alpha_i$ runs from $p$ to some $p_i$
 and $\beta$ is in a ball $B_{q_i}(\delta')$.
 
 Each $g_i$ has a loop $\beta_i$ which is freely r-homotopic to
 $\alpha_i^{-1}\beta_i\alpha_i$ and $\beta_i$ is contained in a $\delta'$ ball
 but it doesn't necessarily have length less than $2 \delta'$
 However, by Lemma~\ref{balltoL}, each $\beta_i$ is freely r-homotopic
 to a product of curves of length $<2\delta'<2\delta_{k+1}$.
 
 Thus $h$ is a product of elements of $\bar{\pi}_1$ which have representative
 curves freely 
 r-homotopic to curves of length $<2\delta_{k+1}$ so $m$ of these
 elements is $<2\delta_{k+1}$.  
 
 We also know by the definition of $\delta_{k+1}$ (\ref{eqdefdel}) that for 
any
 $g \in \bar{\pi}_1 \setminus G_k$ we have  $m(g) \ge \delta_{k+1}$,
 
 Thus $h$ is a product of elements in $G_k$ and $h$ itself is in $G_k$.
 This is a contradiction.
 \newline
 
 To show (\ref{deltaspec}) we must show that if $\delta\in
 (\delta_k,\delta_{k+1})$ 
 then $\delta$ is not in $CovSpec(X)$.
 
 Suppose $\delta\in (\delta_k, \delta_{k+1})$ is in $CovSpec(X)$.
 By Lemma~\ref{lemdelpair}, there exists $g\in \bar{\pi}_1(X)$ with 
 $f(g)=\delta, m(g)=2\delta$. So $\frac 12 m(g) < \delta_{k+1}$ and
 $g$ must be in $G_k$.  But then $f(g) \le \delta_k$ which is a 
 contradiction.
 \newline
 
 Now let $G_{k+1}$ be the group generated by $G_k$ and elements
 $h\in \bar{\pi}_1$ such that $m(h)=2\delta_{k+1}$. To finish c) we need to
 show $G_{k+1}$ includes all $g$ such that $f(g)=\delta_{k+1}$.
 
 Let $h \in \bar{\pi}_1$ be an element such that $f(h)=\delta_{k+1}$
 Then for any $\delta'> \delta_{k+1}$ 
 $h$ lifts trivially to $\tilde{X}^{\delta'}$. By Lemma~\ref{mfinite} we can 
 choose
 \be
 \delta_{k+1} < \delta'< (1/2) \min( \{m(g): g\in \pi_1\} \cap (2\delta_{k+1},
\infty) )
 \ee
 so that if $m(g) < 2\delta'$ then $m(g) \le 2\delta_{k+1}$.
 
 Now $h$  is a product of elements $g_1g_2 \cdots g_k$ where
 each $g_k$ has a representative curve based at $p$ of the form
 $\alpha_i^{-1}\beta_i \alpha_i$ where $\alpha_i$ runs from $p$ to some $p_i$
 and $\beta$ is in a ball $B_{q_i}(\delta')$.
 
 Each $g_i$ has a loop $\beta_i$ which is freely r-homotopic to
 $\alpha_i^{-1}\beta_i\alpha_i$ and $\beta_i$ is contained in a $\delta'$ ball
 but it doesn't necessarily have length less than $2 \delta'$.
 However, by Lemma~\ref{balltoL}, each $\beta_i$ is freely r-homotopic
 to a product of curves of length $<2\delta'$ and, by the
 choice of $\delta'$, to a product of curves of length $\le 2\delta_{k+1}$.
 
 Thus $h$ is a product of elements of $\bar{\pi}_1$ which have representative 
 curves freely r-homotopic to curves of length $\le 2\delta_{k+1}$ so $m$ of 
 these
 elements is $\le 2 \delta_{k+1}$. 
 
 Thus $h$ is a product of elements in $G_{k+1}$ and $h$ itself is in 
 $G_{k+1}$.
 
 This finishes the proof of the induction hypothesis.
 
 Finally using the finiteness of the covering spectrum [Lemma~\ref{univ}], we
 know
 that this process must terminate.  Thus by (c), eventually $G_k$ must equal
 $\bar{\pi}_1$.  So we've determined the value of $f$ for every element
 of $\pi_1$ and determined the marked covering spectrum of $X$.
 \qed
 
 The following examples demonstrate that the length spectrum
 alone does not determine the covering spectrum.  We have many more examples 
in
Section~\ref{laplace} which have the same Laplace spectra and length 
spectra,
but different covering spectra.  
 
 \begin{ex} \label{samelap}
 Let $M_1= S^2_{\pi/2}$ be the standard sphere of diameter
 $\pi/2$, $M_2= \mathbb RP^2_\pi=S^2_{\pi}/\mathbb Z_2$.  Then the length
spectra of both $M_1$ and $M_2$ are 
  $\{ l\pi:  l\in \mathbb N\}$, 
while the covering spectrum of $M_1$
 is empty and the covering spectrum of $M_2$ is $\{\pi/4\}$.
Here $M_1, M_2$ have different fundamental groups. 

There are also examples 
with
same fundamental group.
Let  $M_1=  \mathbb RP^2_\pi \times S^2_{\pi}$
and $M_2=  \mathbb RP^2_{2\pi} \times S^2_{\pi/2}$. 
Their length spectrum is 
 \be
\{\sqrt{(k\pi)^2+(2l\pi)^2} :  k,l \in \mathbb N \cup \{0\}\}\setminus \{0\}.
\ee
 The covering spectrum of $M_1$ is $\{\pi/4\}$ while the covering spectrum 
of
$M_2$ is $\{\pi/2\}$.
 \end{ex}

 %*************************************************************
 
 \sect{Counting generators of fundamental groups}
\label{sectmult} 

 In this section we restrict ourselves to compact length spaces $X$ which have
 universal covers.
 
 The sequence of groups $G_i$ and sets $S_i$ in Theorem~\ref{mlscs} give us a 
way to construct 
 a short basis of $\bar{\pi}_1(X)$ and to define the multiplicity of the 
covering spectrum.

 \begin{defn}  \label{multi}  {\em
  For each $\delta_j \in CovSpec (X)$, the {\em basis multiplicity} of 
 $\delta_j$ is the minimum number of $g\in S_j$ required to generate $G_j$.
 
 Let $\bar{S}_j\subset S_j$ be a list of such generators.  Let
 a {\em short basis} of $\bar{\pi}_1(X)$ be $S=\bigcup \bar{S}_j$.
 }
\end{defn}
 
 Note that the covering space of a compact length space with lifted metric is 
 a locally
 compact complete length space, therefore by the Hopf-Rinow theorem for metric
 space (see \cite{Gr} or \cite[Theorem 2.5.28]{BBI}) each bounded closed 
domain
 is compact. Hence we show below that the multiplicity in above is always 
finite
for 
 compact length space. In fact we have the following lemma:
 
 \begin{lem} \label{tilNi(a,b)}
 Let $X$ be a compact length space with a universal cover $\tilde{X}$.
 For $b \ge a > 0$ and $D=diam(X)$, 
 \be \label{covinab}
 \#_m\{CovSpec(X)\cap [a,b]\} =\sum_{j: \delta_j \in [a,b]} \#\bar{S_j}
  \le \tilde{N}(a,2b+2D+a),
 \ee
 where $\tilde{N}(a,b)$ stands for the number of maximal disjoint balls 
 of radius $a$ fit in a ball of radius $b$ in $\tilde{X}$.
 \end{lem}
 
 This estimate in particular gives an estimate on the multiplicity
of a fixed element $\delta$.
 Lemma~\ref{tilNi(a,b)} will be improved later, see Corollary~\ref{improved}.
 
 \noindent {\bf Proof of Lemma~\ref{tilNi(a,b)}}. Let $\{\lambda_1,\cdots,
 \lambda_k\} = CovSpec(X)\cap [a,b]$ counted with multiplicity. By
Lemma~\ref{lemdelpair} for 
 each $\lambda_i$, there is a $g_i$ in $\bar{\pi}_1(X)$ such that $a\le 
f(g_i)=
 \frac 12 m(g_i) \le b$. By the proof of Theorem~\ref{mlscs}, $\frac 12
 m(g_i^{-1}g_j) \ge f(g_i^{-1}g_j) \ge a$ for all $i\not= j$. Fix $\tilde{p} 
 \in
 \tilde{X}$, we have $d(\tilde{p}, g\tilde{p})) \le m(g)+2D$ for any $g \in
 \bar{\pi}_1(X)$. Therefore each ball $B(g_i\tilde{p},a)$ is disjoint from 
each
 other for $i= 1, \cdots,k$ and all are isometric
 and lie in the ball $B(\tilde{p}, 2b+2D+a)$. This gives (\ref{covinab}).
 \qed

The following example shows that the multiplicities of short elements
 of the covering spectrum can grow to infinity, while elements
 in the covering spectrum converge to 0.
 
 \begin{ex} \label{handles}
 Let $M^2_j$ be a handlebody with $j$ handles which looks like
 a standard 2 sphere with many small handles on the scale of $1/j$. The
 multiplicity of $1/j$ goes to infinity as $j$ goes to infinity. 
 \end{ex}

 Note that the multiplicity in Definition~\ref{multi} does not agree with the
 multiplicity
 of the length spectrum.  We have deliberately related it to the revised
 fundamental group rather than to free homotopy classes of loops.  This way
 Theorem~\ref{mlscs} 
 immediately gives us the following proposition.
 
 \begin{prop}  \label{mcpg}
 For a compact length space $X$ with a universal cover, $\bar{\pi}_1(X)$ can 
be
 generated by the short basis $S$ of Definition~\ref{multi} and
 $ \# S =\#_m\{CovSpec(X)\}.$
 \end{prop}
 
 Note that the number of generators of a fundamental group
 $\pi_1(X,p)$ may not be finite for a
 compact length space, $X$, with a nonsimply connected universal cover.
 The double cone over the Hawaiian earring is its own
 universal cover, so $\#\{CovSpec X\}=0$, but its fundamental 
 group is uncountable and, in particular, not finitely generated.

%******************************************************************************
 
 \sect{Gromov-Hausdorff convergence and $\delta$ covers}
 \label{sectgh1}

 Here we first prove a convergence
property of $\delta$-covering spaces which
 doesn't hold for universal covers. Then we show that, unlike the length
spectrum,  the covering spectrum
 behaves nicely under Gromov-Hausdorff convergence. 
We begin with the definition of the Gromov-Hausdorff distance
between compact length spaces.

\begin{defn} \label{grdef} {\em \cite[Defn 3.4]{Gr}
Given two metric spaces $X$ and $Y$ the {\em Gromov Hausdorff distance}
between them is defined,
\be
d_{GH}(X,Y)=\inf \left\{ d^Z_H(f(X),f(Y)): \stackrel
{ \textrm{ for all metric spaces }Z, \textrm{ and isometric } }
{ \textrm{embeddings }f: X \rightarrow Z, g: Y \rightarrow Z }
\right\},
\ee
where, $d_H^Z$ is the Hausdorff distance between subsets of $Z$,
\be
d_H^Z(A,B)=\inf\{\epsilon>0: B \subset T_\epsilon(A) \textrm{ and }
A \subset T_\epsilon(B) \}.
\ee
Here $T_\epsilon(A)=\{x\in Z: d_Z(x,A) <\epsilon\}$.

%If $x\in X$ and $y\in Y$, the {\em pointed Gromov-Hausdorff distance}:
%\be
%d_{GH}((X,x), (Y,y))=\inf \left\{ d^Z_H(f(X),f(Y)): \stackrel
%{ \textrm{ for all metric spaces } Z, \textrm{ and isometric embed-} }
%{ \textrm{dings }f: X \rightarrow Z, g: Y \rightarrow Z
%\textrm{ s.t.} f(x)=g(y) }
%\right\}.
%\ee 
}
\end{defn}
  
It is then clear what we mean by the Gromov-Hausdorff convergence of
compact metric spaces.  However, for non-compact metric spaces,
the following looser definition of convergence was defined by
Gromov.

\begin{defn} \label{Gr2def} {\em \cite[Defn 3.14]{Gr}
We say that non-compact length spaces $(X_n, x_n)$ converge in
the pointed Gromov-Hausdorff sense to $(Y,y)$ if
for any $R>0$ there exists a sequence $\epsilon_n \to 0$ such that
$B_{x_n}(R+\epsilon_n)$ converges to $B_{y}(R)$ in the Gromov-Hausdorff
sense.}
\end{defn}

It is easy to see that neither the topology of a metric space nor
the dimension is conserved under Gromov-Hausdorff convergence.  
Two compact spaces are close in the GH sense if they look almost the same
with ``blurry vision'' so that ``small holes'' cannot be seen.  
A sequence of $1 \times 1/j$ tori collapses to a circle losing 
both dimension and topology.  The 
sequence of handlebodies, $M_j^2$, of Example~\ref{handles}
converges to a standard sphere thus losing topology without collapsing
to a lower dimension.  One also can lose regularity as can be
seen when taking a sequence of one-sheeted hyperboloids 
converging to a cone.

 \begin{prop}  \label{delta-precom}
 If a sequence of compact length spaces $X_i$ converges to a compact length 
 space
 $X$ in the Gromov-Hausdorff topology, then for any $\delta >0$ there is a
 subsequence of $X_i$ such that their $\delta$-covers also converges in the
 pointed Gromov-Hausdorff topology.
 \end{prop}
 
 This answers a question in \cite{SoWei2}.  Compare Proposition 3.1 in there. 
 
  By Thm 3.6 in \cite{SoWei}, the limit of the $\delta$-covers (if it exists) 
 is
 always a cover of $X$, but note that two different subsequences could have
 different limits as the next example shows.

 \begin{ex} \label{subseq} Let $X_i$ be tori of side lengths
 $1$ by $(n+1)/(2n)$ alternating with the tori of length
 $1$ by $(n-1)/(2n)$. Then $X_i$ converges to the $1$ by $1/2$ torus. For
 $\delta=1/2$, we get two limits of the $\delta$-covers:
 one is the cylinder and the other is Euclidean space.
 \end{ex}
 
 In the following examples we demonstrate that universal covers
 may not have any converging subsequences.  Recall that Gromov's
 Precompactness Theorem \cite{Gr} states that a set, 
$S$, of compact length spaces
 is precompact iff there is a uniform upper bound, $N(r,R)$, on the
 number of disjoint $r$ balls contained in a ball of radius $R$,
 $N_X(r,R)$:
\be
 \forall r,R>0 \,\,\exists N(r,R)\in {\mathbb{N}}
\textrm{ s.t. we have } N_X(r,R) \le N(r,R) \qquad \forall X \in S. 
\ee  
 
 \begin{ex}\label{dimtoinf}
 Let $M^j_j$ be a flat $j$ dimensional $1\times (1/j)\times (1/j) \times 
\cdots
 \times (1/j)$
 torus.  Then the Gromov-Hausdorff limit of $M^j_j$ is a circle.
 The universal covers of the $M^j_j$ are Euclidean $j$ dimensional
 spaces, so $N_{M_j}(1, 5) \ge 2j$
 and the $M_j$ do not have a converging subsequence.
 \end{ex}

 Other examples include the sequence of spheres
 with small handles, $M^2_j$ in Example~\ref{handles}, which converges in
 the Gromov-Hausdorff sense to the standard
 two sphere and  a sequence of finite sets of circles joined at a commom point
 which converges to the standard Hawaii ring. In
 both cases the sequences of universal covers
 do not having any converging subsequences.
 
  Gromov proved that
 if  $M_j$ are closed manifolds with 
 Ricci curvature uniformly bounded from below and dimension bounded above
 then by the Bishop Gromov Volume Comparison Theorem, $N_{M_j}(r,R)$
 is uniformly bounded \cite{Gr}.  
Since the universal covers of the $M_j$ share these 
curvature and dimension bounds, they do have converging 
 subsequences.  However, even in this case,
 the limits of universal covers are not
 necessarily  covers of the limit space.   An example is a sequence 
 of flat $1 \times 1/j$ tori which collapse to a circle.  
The limit of the universal
covers is the Euclidean plane which is not a cover of a circle.
 
\vspace{1cm}

 \noindent {\bf Proof of Proposition~\ref{delta-precom}}. It's enough to show
 that the set of $\delta$-covers of $X_i$ is
 precompact by finding a uniform bound on $N_{\tilde{X}_i^\delta}(r,R)$.
 Since $X_i$ converge to a limit space $X$ in the GH sense, they
 also converge in the pointed GH sense, so there exists $x_i \in X_i$
 and $x \in X$ such that $(X_i,x_i)$ converges to $(X,x)$. 
 So we need only prove that for all $\epsilon, R>0$, the number
 of disjoint balls of radius $\epsilon$ centered
 in $B_{\tilde{x_i}}(R)$ is uniformly bounded. In fact we can fix
 $\epsilon<\delta<R$ since bounds for these $\epsilon$ and $R$ will
 control the others.

 Let  $\tilde{x}_i \in \tilde{X}_i^\delta$ be a lift of  $x_i$. 
 Let  FD$_i$ be a (closed)
 fundamental domain of $X_i$ based at $\tilde{x}_i$. 
Let the $\epsilon$ almost adgacent generators
\be
F_{\epsilon, i} =\{ g\in G(X_i, \delta): \, g T_\epsilon(\textrm{FD}_i)
 \cap T_\epsilon(\textrm{FD}_i) \neq \emptyset \},
\ee
and, let the adjacent generators be the set
\be
F_{i} =\{ g\in G(X_i, \delta): \, g (\textrm{FD}_i)
 \cap (\textrm{FD}_i) \neq \emptyset \} \subset F_{\epsilon, i}.
\ee

 Now examine $B(\tilde{x}_i, R) \subset \tilde{X}_i^{\delta}$. By Milnor's
lemma
 \cite[Lemma 2]{Mi}, if $d(\tilde{x}_i, g\tilde{x}_i) < \delta k$ for some
 positive integer $k$, then $g$ can be expressed as a $k$-fold product, $g =
 h_1h_2\cdots h_k$, with $h_1, \cdots, h_k \in F_i$. Let $k = [R/\delta] +1$,
 where $[R/\delta]$ is the integer part of $R/\delta$. 
 Thus the number of
 fundamental domains $g$FD$_i$ intersecting $B(\tilde{x}_i, R)$ is bounded by 
 $(\#F_i)^{[R/\delta] + 1}$.
 
 On the other hand if $N_i(\epsilon, D)$ is the number of maximal disjoint
 $\epsilon$-balls in $X_i$, then if $\epsilon < \delta$,  we claim the maximal
 number of disjoint $\epsilon$-balls centered in each fundamental domain 
 FD$_i$ is bounded by
 $N=N_i(\epsilon,D) \cdot \#F_{\epsilon,i}$.  If not, then let 
 $\tilde{y}_1, \cdots, \tilde{y}_{N+1}$  be the centers of these balls
 and $y_1, \cdots, y_{N+1}$ be their projections to $X_i$.  
 Since the covering 
 map
 $\tilde{X}_i^\delta$ is isometric on $\delta$ balls,
 \be
B_{y_j}(\epsilon)\cap B_{y_k}(\epsilon) =\emptyset
\textrm{ iff } 
B_{\tilde{y}_j}(\epsilon)\cap B_{g\tilde{y}_k}(\epsilon) =\emptyset
\,\, \forall g \in G(X_i, \delta). 
\ee
which is equivalent to checking that
\be \label{satthis}
B_{\tilde{y}_j}(\epsilon)\cap B_{g\tilde{y}_k}(\epsilon) =\emptyset
\,\, \forall g \in F_{i, \epsilon}.
\ee
So we can select $N_i(\epsilon,D)$ disjoint $\epsilon$ balls in
$X_i$ by first choosing $y_1$ and eliminating the at most $\#F_{\epsilon,i}$
$y_k$ that fail to satisfy (\ref{satthis}) for $y_1$, then choosing the
next remaining $y_j$ and eliminating the at most $\#F_{\epsilon,i}$
$y_k$ that fail to satisfy (\ref{satthis}) for that $y_j$, and so
on.  This a contradiction.
 
So the total number of balls of radius $\epsilon$ in $B_{\tilde{x}_i}(R)$
is bounded by 
 $(\#F_i)^{[R/\delta] + 1} \cdot N_i(\epsilon,D) \cdot \#F_{\epsilon,i}$. 
Since $F_i \subset F_{\epsilon, i}$, we need
only bound $\#F_{\epsilon,i}$ uniformly in $i$.

 Note that by Theorem 3.4 in \cite{SoWei}, we have surjective homomorphisms
 $\Phi_i: G(X, \delta/2) \ra G(X_i, \delta)$ for all $i$ large. 
 We can assume that $\diam X_i \le D$.
 If $\bar{\alpha} \in F_{\epsilon,i} \subset G(X_i, \delta)$,
 then it can be represented by a closed curve $\bar{\sigma}$ passing through
 $x_i=\pi(\tilde{x}_i)$ of length $\le 2(D+\epsilon)$.  
From the proof of surjectivity in
  \cite[Theorem 3.4 ]{SoWei}, we can take
 an $\epsilon$ partition of $\bar{\sigma}$ and get a curve
 $\sigma$ passing through $\pi(\tilde{x})\in X$  such that $\Phi_i (\sigma) =
 \bar{\alpha}$ and the length of $\sigma$ is at most
 5 times as long as $\bar{\sigma}$.  
 Thus each element $g \in F_{\epsilon,i} \subset G(X_i, \delta)$ 
 is mapped to by $\Phi_i$ of some
 element $h \in G(X, \delta/2)$ such that 
 $d_{\tilde{X}^{\delta/2}}(h\tilde{x},
 \tilde{x})<10 (D+\epsilon)$.   
 
 Now if $g_1$ and $g_2$ are two distinct elements in $G(X_i, \delta)$ and 
 $\Phi_i
 (h_1) = g_1$ and $\Phi_i (h_2) = g_2$, then $h_1h_2^{-1} \in G(X, \delta/2)$ 
 is
 nontrivial. Any nontrivial element $h \in G(X, \delta/2)$ has
 $d_{\tilde{X}^{\delta/2}}(hx,x) \ge \delta$. 
 So $d_{\tilde{X}^{\delta/2}}(h_1x,h_2x) \ge \delta$.
 
 Hence for all $i$ large
 \be \label{Fibounded}
 \#\{F_{\epsilon,i}\}  \le  \tilde{N}(\delta/2, 10(D+\epsilon)+\delta/2),
 \ee
 where $\tilde{N}_i(\delta/2, R')$ is the maximal number of disjoint
 balls of
 radius $\delta/2$ that fit in a ball,
$B(\tilde{x}, R') \subset
 \tilde{X}^{\delta/2}$, in the limit spaces $\delta$ cover. 
 \qed
 
 An immediate corollary of this is

 \begin{coro}  \label{delta-precompact}
 Let $\cal{M}$ be a GH compact set of length spaces
 and $\cal{M}^\delta$ be the set consisting of their delta covers. 
 Then $\cal{M}^\delta$ is precompact and $\tilde{N}(\epsilon,R)$ is uniformly
 bounded on $\cal{M}^\delta$.
 \end{coro}
 
Note that $\cal{M}^\delta$ need not be compact since a limit
of $\delta$ covers need not be a $\delta$ cover.  See the example
in \cite{SoWei} immediately above Theorem 3.6.

 Corollary~\ref{delta-precompact} 
enables us to give an improvement of Lemma~\ref{tilNi(a,b)}.  
 
 \begin{coro} \label{improved}
 For all $X$ in Gromov-Hausdorff compact set $\cal{M}$ of compact length 
spaces
 with universal covers and $b>a>0$,  $\#_m(CovSpec(X)\cap [a,b])$ is uniformly
 bounded.
 \end{coro}
 
 As another nice application of Corollary~\ref{delta-precompact} we have 
\begin{prop}
The revised fundamental groups of a Gromov-Hausdorff compact set of complete
length spaces with a uniform lower
 bound on their first systole have  finitely many isomorphism classes.
\end{prop}

Here the first systole of $(X) = \inf CovSpec (X)$, which is a natural way of
extending the definition of first systole  to length spaces that aren't
semilocally simply connected.  

This result generalizes Theorem 5 in \cite{ShW} on manifolds. The same proof 
in
\cite{ShW} carries over once we have uniform bound for $\tilde{N}(\epsilon,R)$
on the universal covers.
 \newline
 
 %*********************************************************************
 
\sect{Convergence of the Covering Spectrum} \label{sectgh2}

Note that the length spectrum can change dramatically under Gromov-Hausdorff
convergence as the following examples show.  First we see that
lengths can disappear in the limit.

\begin{ex}  \label{closed-geodesics1}
Let $X_n$ be the boundary of the $1/n$-neighborhood of the closed planar unit
disk in $\mathbb R^3$ with the induced length metric, 
then $X_n$ converges to the
double disk (identification of two closed unit disks along the boundary
circles). The circle $x^2+y^2= (1+1/n)^2$ is a closed geodesic in $X_n$ but 
the
limit curve $x^2+y^2=1$ is not a geodesic in the limit space. In fact
$2\pi(1+1/n) \in $Length Spectrum of $X_n$, but its limit $2\pi$ is not. 
$X_n$
can be easily approximated by $2$ dimensional 
smooth manifolds with same properties. 
\end{ex}

There are also examples where the limit space is a manifold. 

\begin{ex}\label{closed-geodesics2}
Let $Y$ be the flat
$1\times 2$ torus. Let $M_i$ be manifolds constructed from $Y$ as follows:
cut a ball of radius $1/i$ from $Y$ and replaced it with a rescaled piece
$\frac{1}{30i} X$, where $X$ is ball with warped product metric $dr^2+ f(r)^2 
d
\theta^2$, where 
\[ f(r) = \left\{ \begin{array}{ll} sin(r) & r \in [0, \pi/2] \\
                                     1     &  r \in (\pi/2, 28] \\
                                      r-28 &     r \in [29, 30]
\end{array} \right.
\]
 and smooth in between.
Then $M_i$  converge to $Y$ in Gromov-Hausdorff sense, but there are 
$\lambda_i
\in LengthSpec(M_i)$
with $\lambda_i\to 3 \notin LengthSpec(Y)$ since there are smooth closed
geodesics in $M_i$ which converge
to a piecewise geodesic in $Y$ that has two 90 degree angles and wraps
around $Y$ like a ribbon around a gift. 
\end{ex}

Here is an example with the sudden appearance of elements
in the limit's length spectrum far from elements in the sequences
spectra.

\begin{ex}\label{closed-geodesics3}
Let $X_j$ be flat cones capped at both ends with disks 
of radius 1 and 1+1/j respectively, and height 1.  So 
the bigger side get smaller and converges to
a cylinder with both ends capped by disks of radius 1, then $2\pi$ in the
length spectrum of the limit doesn't come from the length spectrum of the
sequence.  This example can easily be smoothed so that both
the sequence of $X_i$ and the limit are diffeomorphic to spheres.
\end{ex}

 Using the result \cite[Theorem 3.6]{SoWei} 
 that the Gromov-Hausdorff limit of the $\delta$-covers of a sequence
 is almost the $\delta$-cover of their limit space, we can show that the 
covering
 spectrum of the sequence and the
 covering spectrum of the limit space is very closely related. Note that the
 counting here is without multiplicity.
 
 \begin{theo} \label{covspecconv}
 If $X_i$ is a sequence of compact length spaces converging to a compact 
length
 space $Y$, then for
 each $\delta \in CovSpec Y$, there is $\delta_i \in CovSpec
 X_i$ such that $\delta_i \ra \delta$. Conversely if $\delta_i \in CovSpec
 X_i$ and $\delta_i \ra \delta>0$, then $\delta \in CovSpec Y$. Moreover,  if 
 the
 universal cover of the sequence and $Y$ exist, then  $\#\{CovSpec X_i\} \geq
 \#\{CovSpec Y\}$ for all $i$ large. 
 \end{theo}
 
  \Pf Let's prove the first statement. If it's not true, there is a $\delta 
\in
 CovSpec Y$ such that no subsequence of CovSpec$(X_i)$ converges to $\delta$,
 namely there exists an $\epsilon>0$  such that CovSpec$(X_i) \cap
 [\delta-\epsilon, \delta +\epsilon] = \emptyset$ for all except finitely many
 $X_i$. So $\tilde{X}_i^{\delta-\epsilon}
 =\tilde{X}_i^{\delta+\epsilon}$ for all except finite many $i$.
 Now, by Proposition~\ref{delta-precom}, a subsequence of the 
 covers converges (for both $\delta-\epsilon$ and for $\delta + \epsilon$).
 Therefore their limits
 $Y^{\delta-\epsilon}=Y^{\delta+\epsilon}$. By \cite[Theorem
 3.6]{SoWei} $\tilde{Y}^\delta =
 \tilde{Y}^{\delta+\epsilon/2}$, contradicting to $\delta \in
 CovSpec Y$.
 
 To prove the second statement, note that $\tilde{X}_i^{\delta_i} \to
 \tilde{X_i}^{\delta'}$
 is nontrivial for all $\delta'>\delta_i$ and
 $\delta_i$ converges to $\delta>0$.  So for all $\delta'>\delta>0$ and 
 $\epsilon\in (0, \delta)$
 we have $\delta-\epsilon<\delta_i<\delta'$  for $i$ sufficiently large and
 $\tilde{X}_i^{\delta-\epsilon} \to \tilde{X}_i^{\delta'}$ is nontrivial.
 Now take the limit as $i \to \infty$ and we get
 $Y^{\delta-\epsilon} \to Y^{\delta'}$ is nontrivial.
 This is true for all $\epsilon\in (0,\delta)$ and $\delta'>\delta$.
 Now by the properties of limit covers \cite[Theorem
 3.6]{SoWei} we have for all
 $\epsilon\in (0,\delta)$ and $\delta'>\delta$,
 $\tilde{Y}^{\delta-\epsilon} \to \tilde{Y}^{\delta'}$ is nontrivial.
 So $CovSpec(Y)\cap [\delta-\epsilon, \delta']$ is nonempty.
 But $CovSpec(Y)$ is discrete at $\delta$, so this forces
 $CovSpec(Y)$ to include $\delta$.
 \qed
 
 An immediate corollary of this is

 \begin{coro} \label{corcovspconv}
 If $X_i$ is a sequence of compact length spaces converging to a compact 
 length
 space $Y$, then the covering spectra
 converge in the Hasudorff sense as subsets of $\mathbb R$:
 \be
 \lim_{i\to \infty} d_H(CovSpec(X_i)\cup \{0\}, CovSpec{Y}\cup \{0\} ) = 0
 \ee
 \end{coro}
 
 \Pf 
 By the definition of Hausdorff convergence (see inside Defn~\ref{grdef}),
 we need only show that for all
 $\epsilon>0$, there exists $N_\epsilon \in \mathbb N$ such that for all $i 
\ge
N_\epsilon$, 
 \ba
 CovSpec(X_i) & \subset & T_\epsilon (CovSpec(Y)\cup \{0\})  \label{XT} \\
CovSpec(Y)  & \subset & T_\epsilon (CovSpec(X_i) \cup \{0\}) \label{TX}
 \ea
 If (\ref{XT}) is not true, then there is an $\epsilon>0$ and a subsequence of
the $i$ such that there exists 
 \be
 \lambda_i \in CovSpec(X_i) \setminus T_\epsilon (CovSpec(Y)\cup \{0\}).
 \ee
 Since the $X_i$ converge to $Y$ they have a uniform upper bound on 
 diameter, $D$, and the $\lambda_i\in (0, D]$, so a subsequence converges
 to some 
 \be
 \lambda \in [0, D] \setminus T_\epsilon (CovSpec(Y)\cup \{0\}).
 \ee
 Thus $\lambda \notin CovSpec(Y)\cup \{0\}$ contradicting
 Theorem~\ref{covspecconv}.
 
If (\ref{TX}) is not true, then there is an $\epsilon>0$ and a subsequence of
the $i$ such that there exists 
\be
\lambda_i \in CovSpec(Y)\setminus T_\epsilon (CovSpec(X_i)\cup \{0\}).
\ee
Since  $\lambda_i\in [\epsilon, D]$, where $D=diam(Y)$, 
and $(CovSpec(Y)) \cap  [\epsilon, D]$ is closed by 
Proposition~\ref{discrete},
a subsequence of $\lambda_i$ converges
to some $\lambda \in CovSpec(Y) \cap [\epsilon, D]$. In particular, for $i$
sufficiently large $|\lambda_i-\lambda|<\epsilon/2$ and
\be
\lambda \in CovSpec(Y) \setminus 
T_{\epsilon/2} (CovSpec(X_i)\cup \{0\}).
\ee
Then by Theorem~\ref{covspecconv}, we know there exists 
$\delta_i\in CovSpec(X_i)$ converging to $\lambda$ which is a contradiction.
 \qed

 Applying Theorem~\ref{covspecconv} to manifolds with Ricci curvature lower 
 bound
 and combining Theorem 1.1 in \cite{SoWei}, we have

 \begin{coro} \label{Ricci-number}
 If $M^n_i$ is a sequence of manifolds with $\Ric \geq (n-1)H$ converges to a
 compact length space $Y$
  then  $\#\{CovSpec M_i\} \geq
 \#\{CovSpec Y\}$ for all $i$ large. 
 \end{coro}
 
 Another useful application of Theorem~\ref{covspecconv} concerns the covering
 spectra of classes of isolengthspectral manifolds:
 
 \begin{theo} \label{connlencov}
 If $\cal{M}$ is a Gromov-Hausdorff connected class of compact length spaces 
 with a common discrete length spectrum, then all compact length
 spaces in $Cl(\cal{M})$
 have the same covering spectrum as well.
 \end{theo}
 
 \Pf
 We need only show that all the spaces in $\cal{M}$ have the same covering
 spectrum.
 Then the same holds true for all compact $Y$ in the closure by
 Theorem~\ref{covspecconv}
 since there will be $X_i$ with uniform covering spectra converging to $Y$.
 
 Suppose there are at least two distinct covering spectra $C_1$ and $C_2$ for
 spaces
 in $\cal{M}$.  Let ${\cal M}_i$ be the subset of $\cal{M}$ of spaces with 
 covering spectra $C_i$.  Clearly these are disjoint sets.  Each is closed as 
a
 subset of $\cal{M}$ by Theorem~\ref{covspecconv}.  Thus we need only show 
each
is relatively open to get a contradiction.
 
Suppose ${\cal M}_1$ is not relatively open.  Then there is a space $Y\in
 {\cal M}_1$ which can be approximated by $X_j\in \cal{M}$ such that
 $CovSpec(X_j)\neq C_1$.  

Thus for each $j$ either there exists $c_j\in C_1\setminus CovSpec{X_j}$
or there exists $c_j\in CovSpec{X_j}\setminus C_1$.

Since $Y$ is compact, there
is a $D$ such that $\diam X_j \le D$.  Furthermore,
all the spaces share the same discrete length spectrum, $L$, 
and length spectra are closed sets.  Thus 
$\frac 12 L \cap (0,D]$ is finite and, by Theorem~\ref{cov-length-spectrum},
$c_j \subset \frac 12 L \cap (0,D]$.  Thus, by the pigeonhole principle,
there exists $c>0$ and a subsequence of the $j$ such that $c_j=c$.

If $c\in C_1\setminus CovSpec(X_j)$ for this subsequence, 
then by Theorem~\ref{covspecconv}, 
there exists $\delta_j\in CovSpec(X_j)$
such that $\delta_j$ converges to $c$.  But these 
$\delta_j \subset (1/2) L \cap [0,D]$, so eventually they must repeat and
we have a contradiction.

Thus $c\in CovSpec(X_j)\setminus C_1$ for this subsequence.  By 
Theorem~\ref{covspecconv} again, 
$c\in CovSpec(Y)=C_1$ (since $c>0$) which is also a contradiction.
 \qed
 
 This leads immediately to the following corollary.

 \begin{coro} \label{corconnlencov}
 If $M_t$ is a one parameter family of compact Riemannian
 manifolds with a common discrete length spectrum (not counting
 multiplicity), then they have the same covering spectrum.
 \end{coro}

%******************************************************

 \sect{Gaps in the covering spectrum}
 \label{sectgap}

 In this section we discuss gap and clumping  phenomenon in the covering
 spectra of compact length spaces. 
 
 Theorem~\ref{covspecconv} immediately gives us the following gap phenomenon 
 near
 0.
 \begin{prop}
 Given a sequence of compact length spaces $X_i$ which converges to a compact
 length $X$ that has a universal cover, there is 
 $\lambda_X>0$ such that for all $0<\epsilon < \lambda_X$, $\exists N_\epsilon
 \in \mathbb N$ such that
 the  Covering Spectum of $X_i$
   has a gap at
  $(\epsilon, \lambda_X)$ for all $i\ge N_\epsilon$: 
 \be
 CovSpec (X_i) \cap  (\epsilon, \lambda_X) = \emptyset. 
 \ee
 \end{prop}
 
 Note that the gap here depends on the limit space. The simplest
 example which illustrates the restrictions on this gap, is a
 sequence of tori collapsing to a circle.  The size of the
 limit circle determines $\lambda_X$ and the speed of collapse
 determines the relationship between $\epsilon$ and $N_\epsilon$.

 In the following we show there are many gaps in the covering spectrum which
 are uniform in size for a compact class of length spaces. Note that a
Gromov-Hausdorff compact set of compact length spaces have a uniform upper
diameter bound.
 
 \begin{prop} \label{gap}
  Let $\cal{M}$ be a Gromov-Hausdorff compact set of compact length spaces 
with
 universal covers and 
 diam $\le D$, and let $S\subset [L_1,L_2] \subset [0,D]$ be a discrete set
which includes the end points $L_1$ and $L_2$, then if 
 \be
 gap_N(X, S)= \mbox{Nth largest element in}\  \{ \lambda_i-\lambda_{i-1} \}
 \ee
among all $ \lambda_i \in (CovSpec(X)\cap[L_1,L_2])\cup S$ in increasing 
order, 
 then $\gap_{\#S-1}(X,S)$ has a uniform lower bound for all $X\in \cal{M}$.
 This lower bound depends on $S$.
 In particular
 \be
 gap_1(X, S)= max \{ \lambda_i-\lambda_{i-1}: 
                    \lambda_i \in (CovSpec(X)\cap [L_1,L_2]) \cup S \}
 \ee
 is uniformly bounded below depending on $S$.
 \end{prop}

 Note that the importance of this result is that the length of the gap 
 interval of the covering spectrum is uniform for all $X \in {\cal{M}}$. 
 On the other hand the exact location of the gap
 can't be uniform as one can see if we take $\cal{M}$ to be the set
 including all flat 2 dimensional tori, circles and the one point space.
 
 Note that if $S=\{0,D\}$ then for simply connected length spaces,
 $gap_{\#S-1}(X,S)=D$, for $X$ with $CovSpec(X)=\{\lambda_1\}$, then 
 $gap_{\#S-1}(X,S)=\max \{\lambda_1, D-\lambda_1\}\ge D/2$.  So its not
 a strong bound for these length spaces.  But then as we progress to length
spaces with large numbers of elements in the covering spectrum, this will 
force
at least one gap which will be significantly larger than
 the average distance between elements.
 
 By taking $S=\{0, D/N, 2D/N... D\}$, we only start getting interesting 
 controls over spaces with more than $N$ elements in the covering spectrum.
 
 Proposition~\ref{gap} implies that there are sequences of gaps
 approaching $0$.  That is, for any $L>0$, there
 exists a $\delta_{{\cal{M}},L}>0$ such that for any $X \in {\cal{M}}$
 $CovSpec(X)$ has a gap of size $\delta_{{\cal{M}},L}$ between $0$ and $L$.
 
 We now prove the gap theorem.  Note that when $S = \{L_1, L_2\}$ with 
$L_1>0$, 
 the lower bound for gap$_1(X,S)$ also follows from Corollary~\ref{improved}. 

 \noindent
 {\bf Proof of Proposition~\ref{gap}:}
 We already know that when $X_i$ converge to $X$ in the GH sense then the
 CovSpec$(X_i)\cup \{0\}$ converges to CovSpec${X}\cup \{0\}$ in the Hausdorff
 sense [Corollary~\ref{corcovspconv}].  So $(CovSpec(X_i)\cap[L_1,L_2])\cup
\{L_1,L_2\}$ converges to
 $(CovSpec(X)\cap[L_1,L_2])\cup \{L_1,L_2\}$ for any $[L_1, L_2]\subset 
[0,D]$.
 
 Since $S$ includes the endpoints, $L_1$ and $L_2$,
  $(CovSpec(X_i)\cap[L_1,L_2])\cup S$ converges to 
 $(CovSpec(X)\cap[L_1,L_2]) \cup S$.
 
 When two discrete sets of numbers in $[L_1,L_2]$ are close in the Hausdorff
 sense, then the gaps are close as well.  That is the largest gaps are
 close, and the second largest and so on.  Eventually, many of the gaps
 will be close to $0$ or nonexistent.  %ADD DETAILS?
 
 So $\gap_{N}(X_i,S)$ converges to $\gap_{N}(X,S)$ as long as $N \le \#S$,
 and in fact converges for all $N$ if we set the gap to 0 when there aren't
 enough elements in the set.
 
 On the other hand, all the covspecs in $\cal{M}$ are closed and discrete and 
 so
 is S, so for each $X$,  $\gap_{\#S-1}(X, S) > 0$.  Since a positive 
continuous 
 function defined on a compact set has a uniform positive lower bound,
 we are done.
 \qed
 
 \begin{ex}\label{nounif}
 If we look at the following  compact set of metric spaces:
 $X_i$  Hawaii ring with rings of radius $1/i^2, 2/i^2. ... (i-1)/i^2, 1/i, 
1$,
 and $X$,    a circle of radius 1,
 then this space is compact and all elements have discrete covering
 spectra.  There is no uniform bound on the number of elements in the
 covering spectra.  
 
 The largest gap in CovSpec$(X_i)$ is $\pi(1-1/i)$ and the rest of the gaps 
 are the same size, $\pi/i^2$.
 
 Taking $S=\{0, \pi\}$, our uniform lower bound on the largest gap exists
 and is $\pi(1/2)$.
 
 Taking $S=\{0,\frac{\pi}{2}, \pi\}$ our uniform bound on the largest gap is
 $\pi(1/2)$ and on the second largest gap is $\pi(1/4)$.  But this second
 largest gap just records the fact that the covering spectra are
 below $\pi/2$.
 
 Taking $S_j=\{0, \pi/j\}$ says more, since we know there is a uniform
 lower bound on the largest gap between 0 and $\pi/j$.  But in fact
 this gap is basically above the majority of the spectra for all but
 finitely many of the $X_i$.
 \end{ex}
 
 This theorem can also be used to show sets of complete metric spaces
 are not compact in the Gromov-Hausdorff sense.
 
 \begin{ex} \label{todense}
 Let $X_j$ be a compact length space formed by $2^j$ circles of radii
 $\{1/2^j, 2/2^j, 3/2^j, \cdots, (2^j-1)/2^j, 1\}$ joined at a common
 point.  Then 
 \be
 CovSpec(X_j)=\{\pi/2^j, 2\pi/2^j, 3\pi/2^j, \cdots, (2^j-1)\pi/2^j, \pi\}
 \ee
 and $gap_1(X, \{0, \pi\})=\pi/2^j$ is not uniformly bounded below.
 Sure enough this sequence has no converging subsequence in the 
 Gromov-Hausdorff sense.
 \end{ex}
 
 Applying Proposition~\ref{gap} to the compact class of manifolds 
 with a uniform lower bound on Ricci curvature, we have
 
\begin{coro} \label{corogap}]
 For all $H\in \mathbb R$, $D>0$, $n\in \mathbb N$, $L>0$, there
 exists a $\delta (H, D, L, n)>0$ such that for any 
  compact manifolds $M^n$ with $\diam (M)\le D,\ \Ric(M^n)\ge (n-1)H$
 there exists $\lambda_M<L$ such that
 \be
 CovSpec(M) \cap [\lambda_M,\lambda_M + \delta (H,D,L,n)]=\emptyset.
 \ee
 \end{coro}
 
 In addition to showing the existence of gaps of a certain size, one
 can study the location of elements in the covering spectrum.  We call
 the following theorem a clumping theorem, since it shows that elements
 in the covering spectra have tendencies to clump around certain locations.
 
 \begin{prop} \label{clumping} If ${\cal{M}}$ is a Gromov-Hausdorff compact
 set of compact length spaces with universal covers and 
 diam $\le D$, then for all
 $\epsilon>0$, there exists $N_\epsilon\in \mathbb N$ and subsets
 $S_1, S_2, \cdots, S_{N_\epsilon} \subset [0, D]$ such that 
 $m(S_i) <\epsilon$ and each $S_i$ is a finite set of intervals
 of the form:
 $$
 S_i=[0, \epsilon_i)\cup \bigcup_{j=1}^{N_i}(d^i_j-\epsilon_i, 
 d^i_j+\epsilon_i)
 $$
 and for all $X \in {\cal{M}}, \exists i \in 1, \cdots, N_\epsilon \,\,\,
 s.t. \,\, CovSpec{X}\subset S_i.$
 \end{prop}
 
 \Pf 
 Let $\cal{D}$ be the set of discrete subsets of $[0, D]$
 which include $\{0\}$.
 Let
 $F: \cal {M} \to \cal {D}$ be defined as $F(X)=CovSpec(X) \cup \{0\}$.
 
 By Corollary~\ref{corcovspconv},
 $F$ is a continuous map when the metric on $\cal{D}$ is
 the Hausdorff metric.  
 
 Now the continuous image of compact set is compact,
 so $F({\cal{M}})$ is compact.  In particular any open cover of $F({\cal{M}})$
 has a finite subcover.
 
 Fix $X \in \cal {M}$,  denote  $CovSpec(X)=\{d_1, d_2,\cdots,d_N\}$. Define
$U_X \subset {\cal D} = B_{F(X)}(r_{h,X})$, where 
    $r_{h,X}= h \min\{d_1, d_2-d_1, \cdots, d_N-d_{N-1} \}$.
 
 Note that $W_h=[0,\epsilon) \cup \bigcup_{j=1}^{N} (d_j-\epsilon,
 d_j+\epsilon)$   is an open subset of $[0,D]$, and
     for $h<1/2$ this is a disjoint collection of intervals.
 For fixed $\epsilon>0$ choose $h>0$ very small so that
 $m(W_h)<\epsilon$.  This determines $r_{h,X}$ for each $X$.
 
 Now $U_X$ form an open cover of $F({\cal{M}})$, so there is a finite 
subcover.
 Let $U_{X_1}, \cdots, U_{X_{N_\epsilon}}$ be that finite subcover
 and let $\epsilon_i=r_{h, X_i}$ and
 $d^i_j$ the $j^{th}$ element in $CovSpec(X_i)$.
 Then
 \be
 S_i=B_{F(X_i)}(\epsilon_i)=U_{X_i}.
 \ee
 
 So for every $X \in {\cal{M}}$ there is an $i\in 1, \cdots, N_\epsilon$ such
that $CovSpec(X)\subset U_{X_i} = S_i$.
 \qed
 
 One can easily see that Example~\ref{todense} also fails to satisfy
 this clumping phenomenon.

 \begin{coro} \label{coroclump}
 For all $H\in \mathbb R$, $D>0$, $n\in \mathbb N$  and for any
 $\epsilon>0$, there exists $N=N(\epsilon, H, D, n)\in \mathbb N$ and subsets
 $S_1, S_2, \cdots, S_{N} \subset [0, D]$ depending on $\epsilon, H, D$
 and $n$ such that 
 $m(S_i) <\epsilon$ and each $S_i$ is a finite set of intervals
 of the form:
 $$
 S_i=[0, \epsilon_i)\cup \bigcup_{j=1}^{N_i}(d^i_j-\epsilon_i, 
 d^i_j+\epsilon_i),
 $$
 and such that for any compact manifold $M^n$ with 
 $\diam (M)\le D,\ \Ric(M^n)\ge (n-1)H$
 \be
  \exists i \in 1, \cdots, N_\epsilon \,\,\,s.t. \,\, CovSpec{M^n}\subset S_i.
 \ee
 \end{coro}

 %******************************************************
 
 \sect{The Laplace spectrum} \label{laplace}\label{sectlap}
  
 In this section we discuss the relationship between the Laplace spectrum
 and the covering spectrum of a compact Riemannian manifold.  Recall that the 
 Laplace spectrum is defined as the set of eigenvalues of the Laplace 
 operator.  The elements of the Laplace spectrum are assigned a multiplicity
 equal to the dimension of the corresponding eigenspace.
 
 It was proven by Colin de Verdiere that the Laplace spectrum determines 
 the length spectrum of a generic manifold \cite{CdV}.  A generic manifold is 
 one with a ``bumpy metric" in the sense of Abraham and, given any Riemannian
 manifold, there is a nearby generic manifold which is close in the $C^5$ 
sense
 \cite{Ab}.  In particular, negatively curved manifolds are
 generic in this sense \cite{Ber}. The generic manifolds are known to have
discrete length spectra \cite{Ber}.    Thus, the Laplace 
 spectrum determines the length spectrum on negatively curved manifolds 
 of arbitrary dimension.  
 
 On Riemann surfaces, Huber proved 
 the length and the Laplace spectrums determine each other completely 
\cite{Hu}. 
Eberlein has 
 shown that on two step nilmanifolds, the marked length spectrum determines 
 the Laplace spectrum \cite{Eb}. 
 
 However, there are pairs of Laplace isospectral manifolds 
first constructed by Carolyn Gordon that have different length spectra when 
 one takes multiplicity into account \cite{Go1}.  
 
 The simplest result we can get from the above is
 
 \begin{prop} \label{isofinite}
 If $\cal{M}$ is a set of Laplace isospectral manifolds which are 
 negatively curved, then there are only finitely many distinct covering 
 spectra for the manifolds in this class.  
 \end{prop}
 
 By Proposition~\ref{mcpg}, this implies that there is a uniform bound
 on the number of generators of the fundamental groups of these manifolds.
 However, this last fact was already known, since this class of manifolds
 is known to have only finitely many homeomorphism classes \cite{BPP}.     
 
 Another application is that complete length spaces in the GH closure
 of $\cal{M}$ have universal covers.  Furthermore $Cl(\cal{M})$ also
 has only finitely many distinct covering spectra and there is a uniform bound
 on the number of generators of the revised fundamental groups of
 these spaces.  
 
 \Pf
 Since $M\in \cal{M}$ are negatively curved and Laplace isospectral, they 
share
 the same
 length spectrum and this length spectrum is closed and discrete.
 They also have a uniform upper bound on diameter by \cite{BPP}.
 The covering spectra are contained in the length spectrum $\cap [0,D]$ by
 Theorem~\ref{cov-length-spectrum}.
 Thus there are only finitely many possible covering spectra.
 \qed
 
 Since, as yet, all known examples of Laplace isospectral sets of manifolds 
 share the same length
 spectrum not counting multiplicity, we make the following conjecture.
 
 \begin{conj} \label{conjlap}
 If $\cal{M}$ is a set of Laplace isospectral manifolds which are 
 with a uniform upper bound on diameter,
 then there are only finitely many distinct covering spectra for the manifolds
 in this class.
 \end{conj}
 
 In the following example we show that the Laplace spectrum does not determine
 the covering spectrum.  In particular we find a pair of Laplace isospectral
 Riemannian Heisenberg manifolds
 which do not share the same covering spectrum.  Note that Pesce has proven 
 that all Laplace 
 isospectral Riemannian Heisenberg manifolds have the same length spectrum not
 counting multiplicities \cite{Ps}.
 
 \begin{ex} \label{isolapnotcov}
 In \cite{Go1}, Gordon studied the Heisenberg 
 manifolds which are of the following form:
 $H_n(\Gamma, g)=(\Gamma \setminus H_n, g)$ where 
 \be
 \Gamma=\Gamma_{\bar{r},\bar{s},c}=\{(\bar{x}, \bar{y},u)\in H_n: \bar{x}\in
 r_1\mathbb Z\times...\times r_n \mathbb Z, 
 \bar{y}\in  s_1 \mathbb Z \times ...\times s_n \mathbb Z, u\in c\mathbb Z \}
 \ee
 where $H_n$ is the (2n+1) dimensional Heisenberg group with multiplication
 \be
 (\bar{x}, \bar{y}, u)(\bar{x}', \bar{y}', u')= (\bar{x}+\bar{x}',
 \bar{y}+\bar{y}', u+u'+\bar{x}\bar{y}')
 \ee
 and the metric $g$ is a diagonal matrix with diagonal $\{a_1, a_2, \cdots, 
 a_n,
 a_1,\cdots,a_n, 1\}$ with
 $0<a_1\le a_2\le...a_n$ at $T_e(H_n)$.
 Note that one needs $\Gamma$ is a subgroup of $H_n$ which is true 
 iff $r_i, s_i \in c\mathbb Z$.
 
 Then the elements of the fundamental group of $H_n(\Gamma,g)$ are elements of
 $H_n$ of the form \\
$(r_1x_1,\cdots, r_nx_n, s_1y_1,\cdots, s_ny_n, cu)$ where $x_i, y_i, u \in
 \mathbb Z$. By \cite{Go1} Cor 2.9, if $x_i$ or $y_i$ is not zero 
 we have the simple formula:
 \be \label{go2.9i}
 m((r_1x_1,\cdots, r_nx_n, s_1y_1,\cdots, s_ny_n, cu))=
 \sqrt{\sum_{i=1}^n a_i(r_i^2x_i^2+s_i^2y_i^2)}.
 \ee
 Otherwise
 \be \label{go2.9ii}
 m((0,\cdots,0,0, \cdots,0,cz))=\min\{|cz|, (4j\pi a_i(|cz|-j\pi a_i) )^{1/2} 
:
j\in
 \mathbb Z, i=1,\cdots,n, 2j\pi a_i<|cz|\}.
 \ee
For a proof of (\ref{go2.9ii}) see \cite{Eb}. 
 Since $(4j\pi a_i(|cz|-j\pi a_i) )^{1/2}$ is increasing in $j$ for
 $0<j<\frac{|cz|}{2\pi a_i}$, we have
 \be \label{m-min}
 m((0,\cdots,0,0,\cdots,0,cz))=\min\{|cz|, (4\pi a_i(|cz|-\pi a_i) )^{1/2},
 i=1,\cdots,n\} 
 \ee
 and 
 \be  \label{mmin}
 m(0,0,cz) \ge m(0,0,c) \ \mbox{ for all integers} \ z. 
 \ee
 
 Note that the elements $(r_ie_i, 0, 0)$ and $(0, s_i e_i, 0)$ generate all 
the
 elements of $\Gamma$
 of the form $(\bar{x}, \bar{y},0)$, so the covering map is determined on 
these
 elements:
 \be
 f((r_1x_1,...r_nx_n, s_1y_1,\cdots, s_ny_n, 0)) \in \{ \frac 12 \sqrt{ 
a_i}r_i,  
 \frac 12 \sqrt{a_i}s_i: i=1,\cdots, n \},
 \ee
 Note also that $(r_ie_i, 0, 0)(0, s_ie_i, 0)=(r_ie_i, s_ie_i, r_is_i)$, so 
 these
 elements
 also generate elements in the center of the form $(0,0, kr_is_i)$.
 
    If $c\not=kr_is_i$ for all $i$ and integer $k$ and  $m(0,0,c) \not\in \{
 \sqrt{a_i}r_i, \sqrt{a_i}s_i: i=1, \cdots, n\}$, 
 by (\ref{mmin}) $f(0,0,cz)=m(0,0,c)$ and the covering spectrum is 
 \be
 \{ \frac 12 \sqrt{a_i}r_i,  \frac 12 \sqrt{a_i}s_i,  \frac 12 m(0,0,c): i=1,
\cdots, n \}.
 \ee
 When there is an $i$ and an integer $k$ such that $c=kr_is_i$ and 
 $m(0,0,c) \ge \max \{\sqrt{a_i}r_i, \sqrt{a_i}s_i\}$ for that particular $i$, 
 then $f(0,0,cz)=\max \{\sqrt{a_i}r_i, \sqrt{a_i}s_i\}$ 
 for that particular $i$ and then the covering spectrum is only $\{ \frac 12
\sqrt{a_i}r_i,   \frac 12 \sqrt{a_i}s_i: i=1, \cdots, n \}$.
 
 This is particularly interesting because Gordon states that
 two Heisenberg manifolds $H_n(\Gamma',g')$ and $H_n(\Gamma, g)$
 are Laplace isospectral iff $ a_i=a'_i$,  $c=c'$ and 
 $\{a_1 r_1 ^2, \cdots, a_n r_n ^2, a_1 s_1 ^2, \cdots, a_n s_n ^2\}$ is a
permutation
 of $\{a_1 (r'_1) ^2, \cdots, a_n (r'_n) ^2, a_1 (s'_1) ^2, \cdots, a_n (s'_n)
^2\}$.
 Thus the only way to get
 an isospectral pair with different covering spectra is to have
 one which includes $m(0,0,c)$ and one which does not.
 
 Let $a_1=1/8$, $a_2=1/2$ and $c=1$, by (\ref{m-min}) $m(0,0,1)=
 (\pi/2(1-\pi/4))^{1/2} \sim .9767$.
 
 If we take $r_1=20$, $r_2=1$, $s_1=10$ and $s_2=1$, 
 then $c=r_2s_2$ and
 \be
 m(0,0,c) \ge \max \{\sqrt{a_2}r_2, \sqrt{a_2}s_2\}=1/\sqrt{2}, 
 \ee
 so
 \be
 CovSpec(H_2(\Gamma, c))=
 \{ 20/(2\sqrt{8})=5\sqrt{2}/2, 
    1/(2\sqrt{2})=\sqrt{2}/4, 
    10/(2\sqrt{8})= 5\sqrt{2}/4\}.
 \ee
 
 If we take $r'_1=2$,  $r'_2=10$, $s'_1=10$, and $s'_2=1$ 
 then $c\neq r_is_i$ for any $i$ and
 \be
 CovSpec(H_2(\Gamma', c))= 
 \{ \sqrt{2}/4, 5\sqrt{2}/2, 5\sqrt{2}/4,
     1/2(\pi/2(1-\pi/4))^{1/2}\}.
 \ee
 \end{ex}

 It is easy to see that one can construct quite a number pairs of
 isospectral Heisenberg manifolds with different covering spectra
 in this manner.  Interestingly Gordon's particular pair
 of isospectral Heisenberg manifolds with different length
 spectrum (counting multiplicity) \cite{Go1}[Ex 2.4 a] do
 share the same covering spectra: $\{1/2, 1\}$.  So there 
 are distinct pairs of Laplace isopectral manifolds that share
 the same covering spectrum.
 
 Next one questions what happens to the covering spectra in a continuous
 family of Laplace isospectral manifolds.  Note that since Pesce
 has shown Laplace isospectral Heisenberg manifolds share the same
 discrete length spectrum, by Theorem~\ref{connlencov}, we know
 that a one parameter family of Laplace isospectral Heisenberg manifolds 
 must share the same covering spectrum.  Thus the two manifolds
 constructed in Example~\ref{isolapnotcov} are not joined by such
 a one parameter family.
 
 The most explored method of constructing Laplace isospectral
 pairs of Riemannian manifolds is using Sunada's method.  Such
 isospectral manifolds are called Sunada isospectal pairs:

 \begin{defn}
 Sunada isospectral pairs of manifolds are pairs of manifolds 
 $M_1=M/H_1$ and $M_2=M/H_2$ with $\pi: M \to M/G$ is a finite
 normal covering and $H_i$ are subgroups of $G$ such that
 for any conjugacy class $G_j\subset G$,
 \be
 \#(G_j\cap H_1)=\#(G_j\cap H_2).
 \ee 
 Sunada proved that these spaces are Laplace isospectral
 and length isospectral.
 \end{defn}

 A special case of Sunada isospectral pairs of manifolds
 are the Komatsu examples \cite[Example 3]{Su}. 
 
 \begin{ex} \label{XKomatsu}
 A Sunada isospectral pair of manifolds is a Komatsu pair
 if $H_1$ and $H_2$ are any pair of finite groups of
 the same order with exponenets of the same odd prime $p$.
 Both are identified with a set $S$ and they are embedded 
 into the symmetric group on $S$  using the left actions
 of the $H_i$ on $S$.  
 
 Now two permutations of the symmetric group are conjugate
 iff they have the same cycle decomposition (c.f. \cite{Her}).
 So a conjugacy class $G_i$ corresponds to a partition
 $\#S=p_1+...+p_k$, where $p_i \ge 1$.  Since $H_i$ have exponents of order
 $p$, they only contain $p$ cycles.  And since they act on the
 left on $S$, their nontrivial elements must move every point
 in $S$, and thus they are complete sets of $p$ cycles
 and they are all in the same conjugacy class:
 $G_1$ corresponding to $\#S=p+p+..+p$.
 So
 \be
 \#(H_1\cap G_1)=\#S-1=\#(H_2\cap G_1)
 \ee
 and $\#(H_i\cap G_j)=0$ otherwise.
 
 Since every symmetric group can be shown to act by isometries
 on some Riemannian manifold, this creates a Sunada isospectral pair.
 In particular, they can be constructed as a Sunada isospectral pair
 whose common finite cover, $M$, is a simply connected compact manifold.
 \end{ex}
 
 \begin{prop} \label{PKomatsu}
 Komatsu pairs of Sunada isospectral manifolds
 share the same covering spectrum 
 which in fact consists of a single element.
 \end{prop}
 
 \Pf
 Let $M_1=M/H_1$ and $M_2/H_2$ be the Komatsu pair with a common simply
 connected finite cover $M$.  Let $M_0=M/G$ where $G$ is the symmetric
 group.  
 
 Now let $m_i: H_i \to \mathbb R$ be the minimum marked length map for
 $M_i$ and $m:G\to \mathbb R$ be the minimum marked length map for $M/G$.
 Note that $m_i(h)=\inf_{x\in M} d_M(x, hx)=m(h)$.  Furthermore
 $m(g_1)=m(g_2)$ whenever $g_1$ and $g_2$ are conjugate because this
 is the minimum length of a loop freely homotopic to a loop representing
 $g_i$.
 
 However, every nontrivial element in either of the $H_i$ is a member of
 the same conjugacy class corresponding to $\#S=p+p+..+p$.  So
 $m_1(h_1)=m(h_1)=m(h_2)=m_2(h_2)$ for all nontrivial $h_i\in H_i$.  Thus the
 covering
 maps are equal as well, and the only element in the covering spectrum
 is this $m(h_i)$.
  \qed
 
 Note that Komatsu pairs of Sunada isospectral manifolds do not
 necessarily have the same covering spectrum counting multiplicity.
 In \cite[Ex 3]{Su}, $H_1=(\mathbb Z/p\mathbb Z)^3$ has three generators and
 thus the only element in its covering spectrum must have multiplicity $3$
 while 
 $$
 H_2=\lp\,  a,b \,\,| \,\,  a^p=b^p=(aba^{-1}b^{-1})^p=e,\, 
 a(aba^{-1}b^{-1})=(aba^{-1}b^{-1})a,\, b(aba^{-1}b^{-1})=(aba^{-1}b^{-1})b \,
\rp
 $$
 has two generators and thus 
 the only element in its covering spectrum must have multiplicity $2$.
 
 %SHOULD WE ADD ANYTHING ELSE TO THIS SECTION PERHAPS ABOUT
 %ISOSPECTRAL FAMILIES...?

 %MAYBE ADD other applications of Theorem~\ref{connlencov} to
 %connected domains in isospectral classes of generic manifolds
 %and other cases
 
 %**********************************************************

 Department of Mathematics and Computer Science,
 
 Lehman College, City University of New York,
 
 Bronx, NY 10468
 
 sormani@g230.lehman.cuny.edu
 
 Department of Mathematics, 
 
 University of California, 
 
 Santa Barbara, CA 93106 
 
 wei@math.ucsb.edu


\begin{thebibliography}{SoWei2}
 \bibitem[Ab]{Ab} R. Abraham, 
 {\em Bumpy metrics.} 1970 
 Global Analysis (Proc. Sympos. Pure Math., Vol. XIV, Berkeley, Calif., 1968) 
 pp. 1--3  Amer. Math. Soc., Providence, R.I. 
 
 \bibitem[Ber]{Ber} M. Berger, {\em Geometry of the Spectrum. I}, Proc. of 
 Symp.
 in Pure Math. 27 (1975), 129-152.

\bibitem[Bes]{Bes} A. Besse, Manifolds all of whose geodesics are closed,
Springer-Verlag, New York, 1978.
 
 \bibitem[BPP]{BPP} R. Brooks, P. Perry, P. Petersen,
 {\em Compactness and finiteness theorems for isospectral manifolds}, J. reine
 angew. MAth. 426 (1992), 67-89.
 
 \bibitem[BBI]{BBI} D. Burago, Y. Burago, S. Ivanov, A course
 in Metric Geometry. Graduate Studies in Mathematics Vol. 33,
 AMS, 2001.
 
  \bibitem[ChCo1]{ChCo1} J. Cheeger, T. Colding,
 {\em On the structure of spaces with Ricci curvature bounded
 below I}, J. Diff. Geom. 46 (1997) 406-480.
 
 \bibitem[ChCo2]{ChCo3} J. Cheeger, T. Colding, {\em On the
 structure of spaces with Ricci curvature bounded below III},
  J. Differential Geom. 54 (2000), no. 1,
37--74.
 
 \bibitem[CdV]{CdV} Y. Colin de Verdiere, {\em Spectre du laplacien et 
 longueurs des géodisoésiques périodiques. I, II.} (French) Compositio 
 Math. 27 (1973), 83--106; ibid. 27 (1973), 159--184. 
 
\bibitem[Cr]{Cr} C. Croke, {\em Rigidity for surfaces of nonpositive
curvature}, Comment. Math. Helv. 65 (1990), no. 1, 150--169.

\bibitem[DuGu]{DuGu}
Duistermaat, J. J.; Guillemin, V. W. 
{\em The spectrum of positive elliptic operators and periodic 
bicharacteristics.}
  Invent. Math. 29 (1975), no. 1, 39--79
(see proc sym pure math XXVII for short version)

 
 \bibitem[Eb]{Eb} P. Eberlein, {\em Geometry of 2-step nilpotent groups with a
 left invariant metric}, Ann. Sci. Ecole Norm. Sup. (4) 27 (1994), 611-660. 
 
 \bibitem[Fa]{Fa}  A. Fathi, {\em Le spectre marqué des longueurs des surfaces
sans points conjugués. (French. English summary) [The marked length spectrum 
of
surfaces without conjugate points] }
C. R. Acad. Sci. Paris Sér. I Math. 309 (1989), no. 9, 621--624. 

\bibitem[Fu]{Fu} K. Fukaya,  {\em
Collapsing of Riemannian manifolds and eigenvalues of Laplace operator.}
Invent. Math. 87 (1987), no. 3, 517--547.


 \bibitem[Go1]{Go1} C. Gordon, {\em The laplace spectra versus the length 
 spectra
 of Riemannian manifolds}, Contem. Math. 51 (1986) 63-80.
 
 \bibitem[Go2]{Go2}
 C. Gordon, {\em Riemannian manifolds isospectral on functions 
 but not on $1$-forms. }
 J. Differential Geom. 24 (1986), no. 1, 79--96.
 
 \bibitem[Go3]{Go3}
 C. Gordon, {\em When you can't hear the shape of a manifold.} Math.
Intelligencer 11 (1989), no. 3, 39--47

\bibitem[Grnt]{Grnt}
R. Gornet
{\em A new construction of isospectral Riemannian nilmanifolds with examples.}
Michigan Math. J. 43 (1996), no. 1, 159--188.

 \bibitem[Gr]{Gr}M. Gromov, Metric structures
 for Riemannian and non-Riemannian spaces, PM 152, Birkhauser, 1999.

\bibitem[Her]{Her} Herstein, Topics in Algebra 2nd edition, Wiley, 1975.
 
 \bibitem[Hu]{Hu} H. Huber, {\em Uber das Spektrum des Laplace-Operators auf 
 kompakten Riemannschen Flachen.}  [On the spectrum of the Laplace operator 
 on compact Riemann surfaces] 
 Comment. Math. Helv. 57 (1982), no. 4, 627--647.
 
 \bibitem[Mi]{Mi} J. Milnor, {\em A note on curvature and
 fundamental group}, J. Diff. Geom. 2 (1968) 1-7.

\bibitem[Ot]{Ot} Jean-Pierre Otal,  Le spectre marqui des longueurs des 
surfaces
courbure nigative. (French) [The marked spectrum of the lengths of surfaces
with negative curvature]  Ann. of Math. (2)  131  (1990),  no. 1, 151--162. 

 
 \bibitem[Ps]{Ps} H. Pesce,  {\em
 Une formule de Poisson pour les varietes de Heisenberg.} 
 C. R. Acad. Sci. Paris Ser. I Math. 315 (1992), no 12, 1279-1281.
 
 \bibitem[ShW]{ShW} Z. Shen, G. Wei, {\em On Riemannian manifolds of almost
 nonnegative curvature}, Indiana Univ. Math. Jour. 40 (1991), 551-565.
 
 \bibitem[SoWei1]{SoWei} C. Sormani and G. Wei, {\em Hausdorff Convergence and
 Universal Covers}, Transactions of the American Mathematical Society 353
 (2001) 3585-3602.
 
 \bibitem[SoWei2]{SoWei2} C. Sormani and G. Wei, {\em Universal Covers for
 Hausdorff Limits of Noncompact Spaces}, Transactions of the American
 Mathematical Society 356 (2004) no. 3, 1233-1270.
 
\bibitem[Sp]{Sp} E. Spanier, Algebraic Topology, McGraw-Hill, Inc., 1966.
 
 \bibitem[Su]{Su} T. Sunada, {\em Riemannian coverings and isospectral
 manifolds.} Ann. of Math. (2) 121 (1985), no. 1, 169--186
\end{thebibliography}
 \end{document}